\documentclass[12pt]{article}
\title{The Sine$_\beta$ operator}
\date{}
\author{Benedek Valk\'o and B\'alint Vir\'ag}

\usepackage{units}
\usepackage{amsmath, amsthm, amssymb}
\usepackage{graphicx}
\usepackage{amsmath, enumerate}
\usepackage{color, url}
\usepackage[longnamesfirst]{natbib}
\bibpunct[ ]{(}{)}{,}{a}{}{,}

    \newtheorem{theorem}{Theorem}
    \newtheorem{lemma}[theorem]{Lemma}
    \newtheorem{proposition}[theorem]{Proposition}

\theoremstyle{definition} 
    \newtheorem{definition}[theorem]{Definition}

\newcommand{\eps}{\varepsilon}
\newcommand{\Z}{{\mathbb Z}}
\newcommand{\ZZ}{{\mathbb Z}}

\newcommand{\UU}{{\mathbb U}}
\newcommand{\R}{{\mathbb R}}
\newcommand{\CC}{{\mathbb C}}

\newcommand{\HH}{{\mathbb H} }

\newcommand{\lstar}{{\raise-0.15ex\hbox{$\scriptstyle \ast$}}}

\theoremstyle{remark} 
\newcommand{\Sineb}{\operatorname{Sine}_{\beta}}
\newcommand{\Sine}{\operatorname{Sine}}
\newcommand{\Airyb}{\operatorname{Airy}_{\beta}}
\newcommand{\Bessel}{\operatorname{Bessel}}
\newcommand{\Sch}{\operatorname{Sch}}

\newcommand{\Sineop}{\mathtt{Sine}_{\beta}}
\newcommand{\Besselop}{\mathtt{Bessel}_{\beta,a}}
\newcommand{\Airyop}{\mathtt{Airy}_{\beta}}
\newcommand{\Huaop}{\mathtt{HP}_{\beta,\delta}}
\newcommand{\Schop}{\mathtt{Sch}_{\nu}}
\newcommand{\Dirop}{\mathtt{Dir}}
\newcommand{\MT}{\mathrm{T}}
\newcommand{\MR}{\mathrm{R}}
\newcommand{\MS}{\mathrm{S}}

\definecolor{violet}{rgb}{0.8,0,0.2}

\newcommand{\ed}{\stackrel{d}{=}}
\newcommand{\cB}{{\mathcal B}}
\newcommand{\mat}[4]{\left( \begin{array}{cc}
#1 & #2  \\
#3 & #4  \\
\end{array} \right)}
\newcommand{\bin}[2]{\left (
\begin{array} {c}
#1 \\
#2
\end{array}
\right )}
\newcommand{\la}{\langle}
\newcommand{\ra}{\rangle}

\newcommand{\ac}{{\text{\sc ac}}}
\newcommand{\dom}{\operatorname{dom}}

\newcommand{\tr}{\operatorname{tr}}
\newcommand{\ind}{\mathbf 1}
\newcommand{\tch}{\upsilon} 

\newcommand{\benedek}[1]{\textcolor{red}{\texttt{#1}}}

\newcommand{\bside}{\noindent\textbf{\benedek{Begin side computation.}}
\begin{footnotesize}}

\newcommand{\eside}{\end{footnotesize}
\noindent \textbf{\benedek{End side computation.}}}
\newcommand{\spec}{\operatorname{spec}}

\begin{document}
\maketitle

\begin{abstract}

We show that $\Sineb$, the bulk limit of the Gaussian $\beta$-ensembles is the spectrum of a self-adjoint random differential operator
\[
f\;\;\mapsto\;\; 2 {R_t^{-1}}\mat{0}{-\tfrac{d}{dt}}{\tfrac{d}{dt}}{0} f, \qquad f:[0,1)\to \mathbb R^2,
\]
where $R_t$ is the positive definite matrix representation of hyperbolic Brownian motion with variance $4/\beta$
in logarithmic time. The result connects the Montgomery-Dyson conjecture about the $\Sine_2$ process and the non-trivial zeros of the Riemann zeta function, the Hilbert-P\'olya conjecture and  de Brange's attempt to prove the Riemann hypothesis. We identify the Brownian carousel as the Sturm-Liouville phase function of this operator.

We provide similar operator representations for several other finite dimensional random ensembles and their limits: finite unitary or orthogonal ensembles, Hua-Pickrell ensembles and their limits, hard-edge $\beta$-ensembles, as well as the Schr\"odinger point process. In this more general setting, hyperbolic Brownian motion is replaced by a random walk or Brownian motion on the affine group.

Our approach provides a unified framework to study $\beta$-ensembles that has so far been missing in the literature. In particular, we connect It\^o's classification of affine Brownian motions with the classification of limits of random matrix ensembles.
\end{abstract}
\newpage

\tableofcontents

\section{Introduction}

\subsubsection*{Operator for the $\Sine_2$ process}

A central theme in random matrix theory is the study of point process limits of random matrix spectra.  Classical models have point process limits that can be characterized by their joint intensity functions. Most famous is the bulk limit of the Gaussian unitary ensemble (GUE), the $\Sine_2$ process.   This point process has remarkable connections to the critical zeros of the Riemann zeta-function.  According to the Montgomery-Dyson
conjecture the set $\mathcal{Z}=\{y\in \R: \zeta(1/2+i y)=0\}$ looks like the $\Sine_2$ process. The conjecture states that
$$(\mathcal{Z}-U t)\log t\Rightarrow \Sine_2$$
in law as
$t\to \infty$ and $U$ is a uniform random variable on $[0,1]$. \cite{Montgomery} as well as \cite{RS} have breakthrough results in this direction.

The Hilbert-P\'olya approach to possibly prove the Riemann hypothesis is built on the conjecture that the set $\mathcal{Z}$ can be realized as the spectrum of a self-adjoint operator (see \cite{Montgomery}).
One of the most famous attempts to carry out a proof is due to de \cite{deBranges}. It is based on the theory of Hilbert spaces of entire functions he developed previously (see also \cite{Lagarias}). This approach would produce a self-adjoint differential operator of the form
\begin{align}\label{123}
\tau f(t)=R(t)^{-1} \mat{0}{-\tfrac{d}{dt}}{\tfrac{d}{dt}}{0} f(t)
\end{align}
acting on functions $f$ mapping an interval to $\R^2$  with a spectrum given by $\mathcal{Z}$. Here $R(t)$ is a positive definite matrix valued function.

In view of the Montgomery-Dyson conjecture it is natural to ask the following question.
\[\text{\it Is there a random self-adjoint operator of the form (\ref{123}) whose spectrum is $\Sine_2$ process?} \]

In this paper we construct such a random operator.

The question whether there is a natural self-adjoint operator where the spectrum has the same distribution as the $\Sine_2$ process has been raised by various authors. See e.g.~\cite{KatzSarnak_book} and \cite{KatzSarnak_AMS}. 
\cite{BO01} and  \cite{MNN}  generalize the notion of eigenvalue to construct a random operator-like object that has generalized eigenvalues given by the $\Sine_2$ process. Our approach is the first to produce an honest-to-goodness natural self-adjoint operator.

\subsubsection*{Operators for $\beta$-ensembles}
The derivation of the bulk limit of GUE relies on the fact that the joint eigenvalue distribution of the finite ensemble can be computed explicitly. The joint density is proportional to the squared Vandermonde determinant $\prod_{i<j} |\lambda_i-\lambda_j|^2$ with respect to an i.i.d.~Gaussian reference measure.
Other solvable classical Gaussian models (GOE, GSE) have similar joint eigenvalue densities involving the first and fourth power of the Vandermonde.
By allowing the exponent of that term to be any $\beta>0$, \cite{Dyson62_1} introduced a  one-parameter extension of the eigenvalue distributions,  called the Gaussian $\beta$-ensemble. Beta versions of other ensembles are obtained similarly. The joint density can be identified with the Boltzmann factor of a one-dimensional log-gas, with $\beta$ playing the role of the inverse temperature (see the monograph \cite{ForBook} for an extensive treatment).

The point process limits of the most important $\beta$-ensembles have been identified in recent years. In \cite{RRV} the soft edge scaling limit of the Gaussian and Laguerre $\beta$-ensembles (the $\Airyb$ process) and in \cite{RR} the hard edge scaling limit of the Laguerre ensemble (the $\Bessel_{\beta,a}$ process) have been derived. The bulk limit of the Gaussian $\beta$-ensemble (the $\Sineb$ process) and the circular $\beta$-ensemble (the $C\beta E$ process) have been derived in \cite{BVBV} and \cite{KS}, respectively. Many of these limits have been shown to be universal, see \cite{BEY}, \cite{BEYedge}, and \cite{KRV}.

The soft and hard edge limit point processes are realized as the spectrum of second order self-adjoint random differential operators, but similar characterization has been missing in the bulk case. We will provide such a characterization using random differential operators of the form (\ref{123}).

 \subsubsection*{Operators and carousels}

 We construct a self-adjoint random differential operator  of the form (\ref{123}) which has an a.s.~pure point spectrum distributed as the bulk limit process $\Sineb$. The proof relies on the Brownian carousel, a geometric
construction of the $\Sineb$ process given in \cite{BVBV}.

The hyperbolic carousel, introduced in \cite{BVBV}, is a geometric functional which maps a real sequence $\{\lambda_k: k\in \ZZ\}$ to the objects $(\gamma,\eta_0, \eta_1,f)$ where  $\eta_0, \eta_1$ are boundary points of the hyperbolic plane $\HH$, $\gamma(t)\in \HH, t\in[0,T)$ is a path and $f:[0,T)\to (0,\infty)$. See Section \ref{s:carousel}, Definition \ref{def:hypc}.

In Proposition \ref{prop:car_dir} of Section \ref{s:carousel}  we show that   the hyperbolic carousel $(\gamma,\eta_0, \eta_1,f)$ can be identified with a self-adjoint differential operator using the classical Sturm-Liouville oscillation theory.  More precisely, under some mild conditions on the path $\gamma$, we construct a differential operator of the form (\ref{123}) with a  pure point spectrum that is the same  as the sequence produced by the carousel. The operator is given as (\ref{123}) with
\begin{align}\label{gendirac}
 R=\frac1{2y} X^t X, \qquad X=\mat{1}{-x}{0}{y},
\end{align}
on an appropriately defined domain with boundary conditions.
 Here $x(t)+i y(t), t\in [0,T)$ is the path $\gamma(t)$ in the upper half plane representation of the hyperbolic plane.  The boundary points $\eta_0$ and $\eta_1$  show up in the definition of the domain as boundary conditions.
The inverse operator $\tau^{-1}$ is a Hilbert-Schmidt integral operator with a finite $L^2$ norm on the appropriate space, with an integral kernel that is explicitly given in terms of $
\gamma, \eta_0$ and $\eta_1$.

In \cite{BVBV} it was proved that if   $\gamma(t)$ is hyperbolic Brownian motion in logarithmic time, $\eta_0$ is a fixed boundary point of $\HH$, and $\eta_1$ as the limit point $\gamma$, then the carousel produces the $\Sineb$ process. See Theorem \ref{thm:bcar} below 
for the precise statement. The connection between hyperbolic carousels and differential operators  provides the appropriate random differential operator for the $\Sineb$ process.
\begin{theorem}\label{thm:intro1}
Let $x(t)+i y(t), t\in [0,1)$ be hyperbolic Brownian motion with variance $\frac{4}{\beta}$, run in logarithmic time $-\log(1-t)$. Let  $\Sineop$ be the differential operator of the form (\ref{123}) with $R$ given by (\ref{gendirac}). Then the operator $\Sineop$ is self-adjoint on an appropriately defined domain, and its spectrum is given by the $\Sineb$ process.
\end{theorem}
The explicit form of the operator is given in Theorem \ref{thm:SSO} of Section \ref{s:bcar}.

In Section \ref{s:KS} we use the connection between hyperbolic carousels and differential operators to show that the point process scaling limit of the circular $\beta$-ensembles (the $C\beta E$ process) is the same as the $\Sineb$ process. \cite{Nakano} has recently proved this equivalence by deriving both processes as the limit of the same sequence of models.

\subsubsection*{Unitary matrices and Dirac operators}

We prove that Dirac operators of the form \eqref{123} can be associated to finite unitary matrices in a natural way. In Section \ref{s:unitary} we show the following theorem.
\begin{theorem}\label{thm:intro2}
Suppose that $V$ is an $n\times n$ unitary matrix with distinct eigenvalues $e^{i \lambda_k}, 1\le k\le  n$. Then there is a self-adjoint differential operator of the form (\ref{123}) with   spectrum that is equal to the set $\Lambda=\{n \lambda_k+2\pi n j: 1\le k\le n, j\in \Z\}$.  Moreover, there is a hyperbolic carousel driven by a piecewise constant path that produces the same point sequence $\Lambda$.
\end{theorem}
See Propositions \ref{prop:unitary_car}, \ref{prop:unitary_car2}, and \ref{prop:unitary_car2b} in Section \ref{s:unitary} for the exact statements. Our approach builds on the Szeg\H o recursion of the orthogonal polynomials on the unit circle. The path of the constructed is built from the Verblunsky coefficients, the main ingredients of the Szeg\H o recursion.

\cite{KillipNenciu} constructed random unitary matrices with a spectrum that is distributed as the circular $\beta$-ensemble, the $\beta$-generalization of random unitary matrices. Using their results and Theorem \ref{thm:intro2} we  construct random differential operators of the form (\ref{123}) with spectrum distributed as  the circular $\beta$-ensemble. Formal rescaling of these differential operators yield the $\Sineop$ operator of Theorem \ref{thm:intro1}. The rigorous derivation of this operator level convergence is carried out  in the forthcoming paper \cite{BVBV_future}.

We also consider a generalization of the circular $\beta$-ensemble called the Hua-Pickrell ensemble (or circular Jacobi $\beta$-ensemble). Using results of \cite{BNR2009} we construct the corresponding random differential operator. These results are stated in Propositions \ref{prop:circ_car} and \ref{prop:Hua} in Section \ref{s:circ}.

\subsubsection*{Classification of operator limits of random matrices }

The hyperbolic plane $\HH$ can be identified with the affine group of matrices of the form $$X=\mat{1}{-x}{0}{y},$$ which in turn can be identified with the $2\times 2$ positive definite matrices of determinant $1$ via the correspondence $X\mapsto R=\frac{X^t X}{2\det X}$. This transforms a path $x+i y$ in the upper half plane to a path on the affine group of matrices and to a positive definite $2\times 2$ matrix valued function $2R_t$ with $R_t$ given in (\ref{gendirac}).

The random differential operator for the $\Sineb$ process is built from hyperbolic Brownian motion, while the operators corresponding to the finite circular ensembles are built from random walks on the hyperbolic plane.

In several examples we represent point process limits of classical random matrix models (and their $\beta$ generalizations) as the spectrum of similar random differential operators. These operators are built from diffusions on the hyperbolic plane which in turn correspond to right Brownian motions on the affine group of $2\times 2$ matrices.

In particular, we show that the hard edge limit operator constructed in \cite{RR} is equivalent to a differential operator of the form (\ref{123}), and the corresponding hyperbolic carousel is driven by a real Brownian motion with drift embedded in the hyperbolic plane.

Finally, we obtain the
$\Airyb$ process constructed in \cite{RRV} from a  $2\times 2$ canonical system, a generalization of operators of the form (\ref{123}).

\section{Quick review of hyperbolic geometry} \label{s:hyp}

We give an of the hyperbolic plane and the hyperbolic Brownian motion.  For a more detailed discussion see e.g.~\cite{hypgeo} and \cite{hypBM}.

\subsection{The hyperbolic plane}


We will focus on  two models of the two dimensional hyperbolic geometry. The first is the upper half plane $\{z\in \CC: \Im z>0\}$, for which (with a slight abuse of notation) we will also use the notation $\HH$. The second is the
the Poincar\'e disk model  $\UU=\{|z|<1, z\in \CC\}$.

The boundary points of the hyperbolic plane are represented by  $\R\cup \{\infty\}$ for $\HH$ and  $\{|z|=1\}$ for $\UU$. The lines in both models are circular arcs or Euclidean lines that are perpendicular to the boundary. Angles are measured the same way  as Euclidean angles. Distance along a line is measured  by integrating $(\Im z)^{-1}$ or $\frac{2}{1-|z|^2}$, respectively.
For the half plane model the distance between two points $x_1+i y_1$ and $x_2+i y_2$   can be expressed explicitly as
\[
 d_\HH(x_1+i y_1, x_2+i y_2)=\operatorname{arccosh}\left( 1+\frac{(x_1-x_2)^2+(y_1-y_2)^2}{2 y_1 y_2}\right).
\]
The two models can be mapped into each other using the Cayley transform or its inverse. The version we use will map $i\mapsto 0$ and $\infty\mapsto 1$  (some authors use the version  $\infty \to -1$).

The transform is given by the linear fractional transformation
\begin{align}\label{Cayley}
U: \HH\to \UU, \qquad U(z)=\frac{z-i}{z+i}, \qquad U^{-1}(w)=i\frac{w+1}{-w+1}.
\end{align}
Note that the map extends to the boundaries as  $r\to e^{i 2 \rm{arccot}(r)}$ and $e^{i \theta}\to -\cot(\tfrac{\theta}{2})$. In particular, for the hyperbolic angle in $\mathbb H$ with $r\in \partial \mathbb H= \mathbb R\cup{\infty}$ we have
\begin{equation}\label{eq:hangle}
\theta=\mbox{angle}(\infty,i,r) = -2\rm{arccot}(r)
\end{equation}
often, we write $r=a/b$, then $\theta$ is $-2$ times the usual Euclidean angle of the vector $(a,b)$.

For a boundary point $\xi$ and two points $a, b$ in the hyperbolic plane we define the \emph{horocyclic distance} as
\begin{align}\label{def:horo}
d_\xi(a, b)=\lim_{z\to \xi} \left| d(a, z)-d(b,z)\right|.
\end{align}
This is well-defined, and there are explicit formulas in both models. We record the half-plane formulas
with  $\xi=q\in \R$ and $\xi=\infty$:
\begin{align}\label{horo_H}
d_q(x+i y, i)=\log\left(\tfrac{(x-q)^2+y^2}{(1+q^2) y}\right), \qquad
d_\infty(x+i y, i)=\log\left(\tfrac{1}{y}\right).
\end{align}


The orientation preserving isometries of both models can be described by M\"obius transformations  $z\to \frac{az+b}{cz+d}$ with $a,b,c,d\in \CC, ad-bc=1$. In the half plane model we have $a,b,c,d\in \R$ while in the unit disk model we have $c=\bar b$, $d=\bar a$.

Let $
\mathcal P \bin{r_1}{r_2}=\frac{r_1}{r_2}$ be the projection  operator $\mathcal P: \CC^2\to \CC$. If  $M=\mat{a}{b}{c}{d}$ with $\det M\neq 0$ and $\mathcal P \bin{r_1}{r_2}=z$ then
$
\mathcal P M \bin{r_1}{r_2}=\frac{az+b}{cz+d}$.
This way the isometries of the half plane and disk models correspond to elements of SL$(2,\R)$ and SU$(1,1)$ respectively.

The matrix $$\tilde U=\frac1{\sqrt{2}} \mat{-1}{i}{1}{i}$$ corresponds to the Cayley transform $U(z)$ given in (\ref{Cayley}) and can be used to transform elements of $\rm{SU}(1,1)$ to elements of $\rm{SL}(2, \R)$ (and vice versa). Let $B^A=A^{-1}BA$ for the conjugation of matrices. Then for $M\in\rm{SU}(1,1)$  we have $M^{\tilde U}\in \rm{SL}(2,\R)$.

The non-trivial orientation preserving isometries of $\HH$ are classified according to how many fixed points they have. An isometry can have a single fixed point inside the plane, a single fixed point on the boundary or two fixed points on the boundary.

The case with a single fixed point inside the hyperbolic plane is called a hyperbolic rotation.
If the fixed point is 0 in the Poincar\'e disk representation then the hyperbolic rotation this is the same as Euclidean rotation: $z\to e^{i\theta} z$. The general hyperbolic rotation can be obtained by conjugating this `classical' rotation with a M\"obius transformation of $\UU$. The hyperbolic rotations in the half plane model are obtained after by conjugation with the inverse of the Cayley transform.

The evolution of a boundary point under a rate $\lambda$ rotation is as follows. In the disk model when the center of rotation is zero the evolution is described by the ODE $z'(t)=i \lambda z(t)$.  In the half plane model and the center of rotation is $x+i y$  the appropriate conjugation can be applied to show that the evolution is given by
\begin{align}
r'(t)=\lambda \frac{y^2+(r(t)-x)^2}{2y}. \label{hyprot}
\end{align}
For $\lambda>0$ the function $r(t)$ restarts at $-\infty$ whenever it hits $\infty$, for $\lambda<0$ it restarts at $\infty$  whenever it hits $-\infty$.  \medskip

\subsection{The hyperbolic plane, the affine group, and positive definite matrices} \label{subs:hypmat}

The group SL$(2,\mathbb R)$ acts on the upper half plane $\mathbb H$ by M\"obius transformations. This is also true for the group of all $2\times 2$ matrices with positive determinant. An important subgroup of these matrices
are the ones where the corresponding M\"obius transformation  fixes the point $\infty$. This group is the product of the group of constant matrices $cI,\ c>0$ (all of which act
trivially on $\mathbb H$), and the affine group of matrices of the form
\begin{equation}\label{affine1}
X=\mat{1}{-x}{0}{y}, \qquad y>0, x\in \mathbb R.
\end{equation}
The minus sign in $-x$ is unimportant but will be convenient later. In the half-plane action the matrix $X$ takes the point  $x+iy$ to $i$.
Formula \eqref{affine1} gives a correspondence between points $x+iy$ in the upper half plane (or equivalently, the hyperbolic plane) and the affine group.

The matrix
\begin{equation}
R=\frac{X^tX}{\det X}
\end{equation}
is a positive definite matrix of determinant one; the map $X\mapsto R$ is a bijection
between the affine group and all positive definite  $2\times 2$ matrices of determinant one. We call $R$ the {\bf positive definite representation} of the point $x+iy$ in (the Poincar\'e half--plane model of) the hyperbolic plane $\mathbb H$.

\subsection{Hyperbolic Brownian motion}\label{subs:hypB}

In this section we review some basic properties of Brownian motion on the hyperbolic plane,  a conformally invariant diffusion process.  See \cite{hypBM} for proofs and additional details.

\begin{definition}
In the half-plane model standard hyperbolic Brownian motion is the solution of the SDE
\begin{align}
d\mathcal B = \Im \mathcal B\, dZ \label{hypBM01}
\end{align}
where $Z$ is standard complex Brownian motion, i.e. $\Im Z$, $\Re Z$ are independent standard real Brownian motions.
If we replace $dZ$ by $\sigma dZ$ we get hyperbolic Brownian motion with variance $\sigma^2$.
\end{definition}

We note that an equivalent SDE in the Poincar\'e disk model is
$
d\tilde \cB=\tfrac{1}2 (1-|\tilde \cB|^2)\,dZ.
$

With $\mathcal B(0)=i$  the solution  of the
SDE \eqref{hypBM01} is
\begin{align}\label{explicit}
\Im \cB(t)= e^{\Im Z(t)-t/2}, \qquad \Re \cB(t)= \int_0^t e^{\Im Z (s)- s/2} d\Re Z(s).
\end{align}
In the affine representation \eqref{affine1} this reads
$$
dX=\mat{0}{d\Re Z}{0}{d\Im Z} X, \qquad X=\mat{1}{-\Re \cB}{0}{\Im \cB}
$$
which is simply a right Brownian motion on the affine group \eqref{affine1}. The minus sign ensures that corresponding fractional linear transformation takes $\cB \in \mathbb H $ to $i\in \mathbb H$.

Hyperbolic Brownian motion converges to a boundary point $\cB(\infty)$ as $t\to \infty$. With $\cB(0)=i$ we have
\[\cB(\infty) = \int_0^\infty e^{ \Im Z(s)-s/2} \,d\Re Z(s).
\]
The process $\cB$ is invariant under hyperbolic rotations fixing $\cB(0)$, and so $\mathcal B(\infty)$ has a rotationally invariant distribution. We say that $\mathcal B(\infty)$ is
is {\bf uniformly distributed as seen from $\mathcal B(0)$}. By the Markov property the following also holds. For any fixed $T$ the random variable $\cB(\infty)$ conditionally on $\cB(t), t\in [0,T]$ is uniformly distributed as seen from $\cB(T)$.

Hyperbolic Brownian motion conditioned to converge to a given boundary point is also a diffusion, as can be seen from a standard application of Doob's $h$-transform.  In the half-plane model the process conditioned to converge to $\infty$ satisfies the SDE
$$d\cB =\Im \cB (dZ+i dt).$$
Finally, we describe reversed hyperbolic Brownian motion on a finite time interval, as viewed from the endpoint. Let $\cB(t), t\in[0,s]$ be hyperbolic Brownian motion in the half-plane model then the process
\begin{equation}\label{eq:bhat}\hat \cB(t)=\frac{\cB(s-t)-\Re \cB(s)}{\Im \cB(s)}, \qquad t\in[0,s]
\end{equation}
is a time reversal, translated to start at $i$ by the hyperbolic isometry fixing the boundary point $\infty$. Using the explicit solution \eqref{explicit} we see that $\tilde B$ has the same distribution as hyperbolic Brownian motion started from $i$, conditioned to hit $\infty$.

\begin{proposition}[Time-reversal symmetry of hyperbolic Brownian motion]\label{prop:reverse}
Fix a time $s$, and given $\cB(t),t\in[0,s]$, let $X$ be a uniform random point on $\partial \mathbb H=\mathbb R\cup\{\infty\}$ as seen from $\cB_s$. Let $\MT$ be the M\"obius transformation taking $\cB(s)$ to $i$ and $X$ to $\infty$. Let $\tilde \cB(t)=\MT\cB(s-t)$. Then 

\begin{equation}
\label{eq:timereversal}
((\tilde \cB(t)\in[0,s]),\MT\infty)\ed ((\cB(t),t\in[0,s]),X).
\end{equation}
\end{proposition}
\begin{proof}

Consider the process $\hat \cB$ in \eqref{eq:bhat}.  It is a time-reversal of $\cB$ mapped by the M\"obius transformation  $S$ that fixes $\infty$ and moves $\cB(s)$ to $i$. The law of $\hat \cB$ is Brownian motion conditioned to converge to  $\infty$. We extend $\hat \cB$ from $[0,s]$ to $[0,\infty)$ to get a hyperbolic Brownian motion conditioned to converge to $\infty$.

We can write $\MT=\MR\MS$ where $\MR$ is a uniform random rotation about $i$ independent of $\hat \cB(t), t\in [0,s]$ by the angle $(X,\cB(s), \infty)$. Then $\tilde B=\MR \hat B$ is distributed as hyperbolic Brownian motion. In particular the point $\MT\infty=\MR\hat \cB(\infty)$ given $\tilde \cB(t), t\in[0,s]$ is uniform as seen from $\tilde \cB(s)$. This proves \eqref{eq:timereversal}.
\end{proof}

Note the following heuristic computation. By It\^ o's formula applied to $\tilde \cB =\MR\hat \cB$ where $\MR z=(az-1)/(z+a)$ with $a=\MT\infty$ we get an equation for $\tilde \cB$ as $d\tilde \cB=\Im \tilde \cB d\tilde Z$, where
\begin{equation}\label{eq:driftterm}
d\tilde Z(s-t) = -\frac{ \tilde \cB(s-t)-\MT\infty}{\overline {\tilde \cB(s-t)}-\MT\infty}(dZ(t) + idt)=
-\frac{\overline{ \cB(t)}-X}{ {\cB(t)}-X}(dZ(t)+idt).
\end{equation}
This SDE looks into the future so it is not defined within the It\^o theory. The rotation factor depends on $Z$ so there is no contradiction in the appearance of a drift term on the right. Using enlargement of filtration or rough path theory this computation can be made rigorous.

\section{Dirac operators and canonical systems}\label{s:Dirac}

A canonical system is a differential equation of the form
\begin{align}\label{canonical}
J v'(t)=\lambda R(t) v(t), \qquad J=\mat{0}{-1}{1}{0}.
\end{align}
Here $R(t)$ is an integrable positive semidefinite real $2\times 2$ matrix valued  function and the $\mathbb R^2$-valued function $v(t)$ also depends on $\lambda$.   If $R(t)$ is strictly positive definite then $v(t)$ solves the eigenvalue equation of the differential operator  $\tau$ defined as
\begin{align}\label{Diracop}
\tau  v(t)=R^{-1}(t) J  v'(t). 
\end{align}
Differential operators of this form are called Dirac operators. The goal of this section is to review some of the well-known properties of such Dirac operators, our main source is \cite{Weidmann}.
We consider the operator on the  interval $[0,T)$ with a fixed $T>0$, where $T$ can be $\infty$. \medskip

We will assume that the function $R:[0,T)\to \R^{2\times 2}$ satisfies the following conditions:
\begin{enumerate}[(A)]
\item \label{cond1} $R(t)$ is positive definite for all $t\in [0,T)$.

\item   \label{cond2} $R(t)$ is measurable and $\|R\|, \|R^{-1}\|$ are locally bounded on $[0,T)$.

\item  \label{cond3}  There exists a nonzero vector $u_*\in \mathbb R^2$ for which  $\int_0^T u_*^t R(s) u_* ds<\infty$.

\end{enumerate}

In our applications we will mostly consider the following two examples: $R(t)$ is continuous or piecewise constant on $[0,T)$.

Let $\ac=\ac_{[0,T)}$ denote the space of absolutely continuous functions on $[0,T)$, i.e.~$f\in\ac$ if for some function $g$ we have $f(t)=f(0)+\int_0^t g(s)\,ds$ and $\int_0^t |g(s)|\,ds<\infty$ for all $t\in[0,T)$. For any $f\in \ac$ the value of $\tau f$ is defined almost everywhere on $[0,T)$.

Standard theory of differential equations gives the following proposition. 

\begin{proposition}[see Theorem 2.1 in \cite{Weidmann}]\label{prop:ODE}
For any $a\in \CC^2$, $t_0\in [0,T)$  and $\lambda \in \CC$ there is a unique continuous solution of $\tau v=\lambda v$ with  $v(t_0)=a$. Moreover, the solution $v(t,\lambda)$ is analytic in $\lambda$ for any fixed $t\in [0,T)$.
\end{proposition}

\subsection{Self-adjoint extensions}

We consider the Hilbert space $L^2_R=L^2_ R[0,T)$ with the inner product
\begin{align}\label{L2R}
\la v, w\ra_R=\int_0^T  v(t)^t R(t) w(t)  \,dt.
\end{align}
By assumption (\ref{cond3}), the constant function $f(t)=u_*$ is in $L^2_R$.

Assuming that $v$ and $w$ are differentiable and  have compact support in $(0,T)$, partial integration yields $\la v, \tau w\ra_R=\la \tau v, w\ra_R$.
We will now define
a domain on which $\tau$ is self-adjoint with respect to the inner product \eqref{L2R}. Recall
that a linear operator $A$ with domain $\dom A\subset L^2$ is {\bf self-adjoint} with respect to the inner product of $L^2$ if the following hold:
\begin{itemize}
\item for all $v,w\in \dom A$ we have $\la v,Aw\ra = \la Av,w \ra$, and
\item if $v\in L^2$ satisfies that $w\mapsto \la v,Aw\ra$ is continuous on $\dom A$, then $v\in \dom A$.
\end{itemize}
To define the domain, we will  need to specify boundary conditions.
We distinguish two cases:
\begin{enumerate}
\item[a)] There is a vector $\tilde u$ not parallel to $u_*$ so that the constant function $\tilde u$ is in $L^2_R$. Then all constant vectors are in $L^2_R$ and $\|R\|$ is integrable on $[0,T)$. In this case we say that the operator $\tau$ is   {\bf limit circle} at the right boundary $T$.
\item[b)]  $\|R\|$ is not integrable on $[0,T)$. Then $\tau$ is called {\bf limit point} at the right boundary $T$.
\end{enumerate}
Note that from our assumption (\ref{cond2}) it follows that $\|R\|$ is integrable near 0, thus $\tau$ will always be limit circle at the left boundary $0$.

To motivate the choice of the domain, we recall the following theorem:
\begin{theorem}[Weyl's alternative, Theorem 5.6 in \cite{Weidmann}] 
\label{thm:Weyl}
In the limit circle case, for every $\lambda\in \CC$ all solutions of $(\tau- \lambda)v=0$ are $L^2_R$ near the right boundary $T$.
In the limit point case, for every $\lambda\in \CC$ there exists at most one (up to constant factor) $L^2_R$ solution
of $(\tau- \lambda)v=0$.
\end{theorem}

In the limit point case, the $L^2$ condition on a function forces it to have the same behavior near $T$ as the constant $u_*$. Hence we only need to specify the boundary condition on the left. We define, for some nonzero
vector $u_0\in \mathbb R^2$:
\begin{align}\label{domLP}
\textup{dom}_{LP}(\tau)=\{v\in L^2_R\,:\, v\in \ac, \, \tau v\in L^2_R,\,v(0)^t\,J\,u_0=0\}.
\end{align}
For the limit circle case we need boundary conditions at both endpoints. For nonzero vectors $u_0, u_1\in \mathbb R^2$, we set
\begin{align}\label{domLC}
\textup{dom}_{LC}(\tau)=\{v\in L^2_R\,:\, v\in \ac, \, \tau v\in L^2_R,\,v(0)^t\,J\,u_0=0, \,\lim_{s\to T} v(s)^t\,J\,u_1 = 0\}.
\end{align}
%
%
The expression $a^tJ b=a_2 b_1-a_1 b_2$ is the Wronskian of the constant functions $a$ and $b$, and it vanishes
exactly when the vectors $a$ and $b$ are parallel.

By the standard theory we have:
\begin{theorem}[see~Theorem 5.8 in \cite{Weidmann}]
The operator $\tau$ with the above domain is self-adjoint with respect to $L^2_R$.
\end{theorem}

From this point we fix the domain of the operator $\tau$, this means fixing the boundary condition $u_0$ and $u_1$ in the limit circle case, and $u_0$ in the limit point case. In the limit point case we will use the notation $u_1$ for the vector $u_*$ from condition (\ref{cond3}).  This way the constant function $u_1$ is always the unique (up to a constant factor) solution of $\tau u=0$ in the domain of $\tau$.
We will say that $u_0$ and $u_1$  are the boundary conditions for the operator $\tau$.

\subsection{The inverse operator}

In this section we describe the inverse operator $\tau^{-1}$.
Assume that the boundary conditions $u_0$ and $u_1$ are not parallel. This implies that $0$ is not an eigenvalue of $\tau$, otherwise the constant function $u_0$ would be an eigenfunction and in particular it would be in the domain of $\tau$. If $u_0, u_1$ are not parallel then  $u_0^t J u_1\neq 0$, and because only the 1-dimensional subspaces spanned by $u_0, u_1$ are important in the definitions of the domain, we can assume $u_0^t J u_1=1$.

As the following theorem shows, one can find a simple representation for $\tau^{-1}$.
\begin{theorem}[see~Theorem 7.8 in \cite{Weidmann}] Suppose that $\tau$ is a Dirac operator of the form (\ref{Diracop}) satisfying (\ref{cond1}-\ref{cond3}) and
$u_0^t J u_1=1$.
For any  $v\in L^2_R$   the integral
\[
g(x)=\int_0^T \mathcal K(x,y)R(y) v(y) \,dy,\qquad \mathcal K(x,y)=u_0 u_1^t  \ind(x<y)+u_1 u_0^t \ind(x\ge y)
\]
is finite  for $x\in [0,T)$, and the function $g$ satisfies $\tau g = v$. If $g$ is in $L^2_R$ then $g\in \dom(\tau)$.
\end{theorem}

Thus if $\tau^{-1}$ is a bounded operator (i.e.~if 0 is in the resolvent set of $\tau$) then it is the integral operator with kernel
\begin{align}\label{HSkernel}
\mathcal{K}(x,y)R(y)=\left(u_0 u_1^t  \ind(x<y)+u_1 u_0^t \ind(x\ge y)  \right) R(y).
\end{align}
Note that $R(x) (\mathcal K(x,y)R(y))= (\mathcal K(y,x)R(x))^t R(y)$, hence $\tau^{-1}$ is symmetric
in $L^2_R$ as expected.


The Hilbert-Schmidt norm of $\tau^{-1}$ is given by
$$\|\tau^{-1}\|_2^2=\sum_k \|\tau^{-1} \varphi_k\|_R^2$$
where $\varphi_k$ is an orthonormal basis in  $L^2_R$. A straightforward computation gives the result
\begin{align}\label{HSnorm0}
\|\tau^{-1}\|_2^2=\int_0^T\int_0^T \tr(\mathcal K(x,y)R(y)\mathcal K(x,y)^tR(x))\,dy\,dx.
\end{align}
One way to see this is by conjugating the integral operator $\tau^{-1}$ with a positive definite square root of $R(x)$ to get the integral operator $A=R(x)^{1/2}\mathcal K(x,y)R(y)^{1/2}$ which is now symmetric with respect to the Lebesgue measure. The Hilbert-Schmidt norm of $A$ is given by the usual formula  $\int_0^T\int_0^T \tr(A(x,y)A(x,y)^t)\,dx\,dy$ which simplifies to (\ref{HSnorm0}). If this norm is finite, we call the operator $A$ Hilbert-Schmidt.

By symmetry we can write
\begin{align}\label{eq:HSnorm}
\|\tau^{-1}\|_2^2=2  \int_0^T \int_0^x  u_0^t R(y)  u_0 \, u_1^tR(x) u_1 \,dy\,dx.
\end{align}
This leads to the following theorem.
\begin{theorem}[see~Theorem 7.11 in \cite{Weidmann}]
Assume that the integral  in (\ref{eq:HSnorm}) is finite. Then  $\tau^{-1}$ is Hilbert-Schmidt,   $\tau$ has a discrete spectrum and the eigenvalues satisfy $\sum_i \lambda_i^{-2}<\infty$.
\end{theorem}
By taking the second integral in (\ref{eq:HSnorm}) up to $T$ we get the upper bound
\[
\|\tau^{-1}\|_2^2\le 2 \int_0^T  u_0^t R(y)  u_0 dy \, \int_0^T  u_1^t R(x)  u_1 dx.
\]
In the limit circle case this is always finite and thus the previous theorem always applies.\medskip

We can use the inverse operator to approximate $\tau$ with limit circle type operators. Let $0<T_n\uparrow T$ be a positive increasing sequence approximating $T$. For each $n$ we denote by $\tau_n$ the  restriction of the differential operator $\tau$ to the interval $[0,T_n)$. On $[0,T_n)$ the function $\|R\|$ is bounded, thus $\tau_n$ is limit circle at the right endpoint $T_n$. Then $\tau_n$ is self-adjoint if we set  its domain as
\begin{align*}
\textup{dom}(\tau_n)=\{v\in L^2_R\,:\, v\in \ac, \, \tau v\in L^2_R,\,v(0)^t\,J\,u_0=0, \,v(T_n)^t\,J\,u_1=0\}.
\end{align*}
Note that we used the same boundary conditions as for $\tau$. The inverse operator $\tau_n^{-1}$ is the integral operator defined in  (\ref{HSkernel}), but restricted to $[0,T_n)$. We can view this as an integral operator on $[0,T)$ by setting  the kernel  equal to zero outside of  $[0,T_n)^2$, and the sequence of these integral operators converge to $\tau^{-1}$ in the Hilbert-Schmidt sense.
This gives the following result.
\begin{theorem}[Theorem 1 in \cite{StolzWeidmann93}]\label{thm:appr}
Suppose that  $\tau^{-1}$ is Hilbert-Schmidt, and consider an approximating sequence $T_n\uparrow T$ with corresponding  operators $\tau_n$. Then the eigenvalues of $\tau$ are exactly the limits of the eigenvalues of $\tau_n$ as $n\to \infty$, moreover the corresponding eigenprojections converge in norm. In particular, if $a<b$ are not eigenvalues of $\tau$ then the number of eigenvalues in $[a,b]$ for $\tau_n$ converges to the number of eigenvalues of $\tau$ in $[a,b]$ as $n\to \infty$.
\end{theorem}

\subsection{A parametrization of $\tau$}

Consider the operator $\tau$ from (\ref{Diracop}). Since $R(t)$ is positive definite it has a unique representation in the form of
\begin{align}\label{Rrep}
R=\frac{f}{y} \mat{1}{-x}{-x}{x^2+y^2}
\end{align}
where $f>0$, $y>0$ and $x\in \R$. Specifically,  $f=\sqrt{\det R}$, $y=\frac{f}{R_{1,1}}$ and $x=- \frac{R_{1,2}}{R_{1,1}}$.
As a consequence, we get another useful formula
\begin{align}\label{RXX}
R=\frac{f}{\det{X}}X^tX, \qquad X=\mat{1}{-x}{0}{y}.
\end{align}
If $\det R=1$ then $R$ is exactly the positive definite representation of $x+i y$ introduced in Subsection \ref{subs:hypmat}.

We can use the representation (\ref{Rrep}) to parametrize the operator $\tau$. We introduce the notation $\Dirop(x+i y,  u_0, u_1,f)$ to denote the Dirac operator $\tau$ with $R$ given in (\ref{Rrep}) with boundary conditions $u_0, u_1$. Note that $x+i y:[0,T)\to \HH$ is a function in the Poincar\'e upper half plane, $f:[0,T)\to (0,\infty)$ is a positive function and $u_0, u_1$ are non-zero real vectors. Since only the directions of $u_0, u_1$ matter for the definition of $\tau$, we will use $u_i$ and $\mathcal P u_i  \in \partial \HH$ interchangeably. If $f$ is the constant function $\tfrac12$ then we drop it from the notation, i.e.~$\Dirop(x+i y,  u_0, u_1)=\Dirop(x+i y,  u_0, u_1,\tfrac12)$.

The following lemma summarizes some properties of $\Dirop(x+i y,  u_0, u_1,f)$  in terms of the ingredients $x+i y,  u_0, u_1,f$. Recall from (\ref{def:horo}) the definition of the horocyclic distance of two points in the hyperbolic plane corresponding to a boundary point.

\begin{lemma}
Suppose that the function $x+i y:[0,T)\to \HH$ is measurable and locally bounded, and let $f:[0,T)\to (0,\infty)$ be measurable with $f, f^{-1}$ locally bounded. Let $u_0\neq  u_1$ be boundary points of $\HH$ and  let $r_0$ be a point in $\HH$. Then the $R(t)$ function of the operator $\Dirop(x+i y,  u_0, u_1,f)$ (defined via (\ref{Rrep})) satisfies the  conditions (\ref{cond1})-(\ref{cond3}) if
\begin{align}\label{horo_int1}
\int_0^T f(t) e^{d_{u_1}(r_0, x(t)+i y(t))}dt<\infty.
\end{align}
The inverse operator $\Dirop(x+i y, u_0, u_1,f)^{-1}$ is Hilbert-Schmidt if
\begin{align}\label{horo_int2}
\int_0^T \int_0^t f(s) f(t) e^{d_{u_0}(r_0, x(s)+i y(s))+d_{u_1}(r_0, x(t)+I y(t))} ds dt<\infty.
\end{align}
\end{lemma}
\begin{proof}
By our assumptions on $x+i y$ and $f$ the conditions (\ref{cond1}) and (\ref{cond2}) are immediately satisfied.
Recall the representation (\ref{horo_H}) of the horocyclic distance in the half-plane model. A simple computation shows that  for any nonzero $v\in \R^2$ we have
\[
e^{d_{\mathcal P v}(x(t)+i y(t), i)}=\frac{1}{f(t) |v|^2} v^t R(t) v,
\]
 with $R$ given in (\ref{Rrep}). Thus if (\ref{horo_int1}) holds with $r_0=i$ then  the function $R$ satisfies condition (\ref{cond3}) with $u_*=u_1$. We also have  $|d_\xi(r_0, b)-d_\xi(r_1, b)|\le d(r_0, r_1)$, thus if condition (\ref{horo_int1}) holds for some $r_0\in \HH$, then it holds for any $r_0\in \HH$.

The same argument shows that (\ref{horo_int2}) implies
\[
 \int_0^T \int_0^t u_0^t R(s) u_0 \, u_1^t R(t) u_1 ds \, dt<\infty,
\]
which is equivalent to $\Dirop(x+i y, u_0, u_1,f)^{-1}$ being Hilbert-Schmidt.
\end{proof}

We finish this section by recording a simple change of variables transforming the Dirac operator $\tau$ into a self-adjoint operator on $L^2$, using the representation (\ref{RXX}).

A function $u$ is in  $L^2_R[0,T]$ (see \eqref{L2R}) if and only if $ \tilde X u\in L^2[0,T]$ where $\tilde X=\left(\frac{f}{\det X}\right)^{1/2} X$. The conjugated operator
\begin{align}\label{tildetau}
\tilde \tau u =\tilde X  \tau( \tilde X^{-1} u)=(\tilde X^t)^{-1} J (\tilde X^{-1} u)'
\end{align}
is self-adjoint on the domain $\{u: \tilde X^{-1} u\in \textup{dom}_\tau\}\subset L^2[0,T]$ and has the same spectrum as $\tau$. The inverse $\tilde \tau^{-1}$ is an integral operator on $L^2[0,T]$ with kernel
\begin{align*}
\tilde X(x) \mathcal K(x,y) \tilde X(y)^t=\tilde X(x) \left(u_0 u_1^t  \ind(x<y)+u_1 u_0^t \ind(x\ge y)  \right) \tilde X(y)^t,
\end{align*}
and its Hilbert-Schmidt norm is the same as that of $\tau^{-1}$.

\section{Oscillation theory and the hyperbolic carousel}\label{s:carousel}

\subsection{The phase function} \label{subs:phase}

For $\lambda\in \R$ let $v(t,\lambda), t\in[0,T)$ be the solution of
\begin{align}\label{osc}
\tau v=\lambda v \qquad \textup{ with } \quad v(0,\lambda)=u_0.
\end{align}
Then $v\in\mathbb R^2$ is continuous on $[0,T)\times \R$ and by Proposition \ref{prop:ODE} it is never equal to $(0,0)^t$.

Since only the 1-dimensional subspaces spanned by the vectors $u_0$ and $u_1$ are relevant for the definition of the operator, we may assume that $u_i=(\cos(\varphi_i/2), -\sin (\varphi_i/2))^t$ with $\varphi_i\in [0,2\pi)$ for $i=0,1$.   Let $\phi_\lambda(t)$ be twice the angle of $(1,-i)\cdot v(t,\lambda)$, that is the unique real-valued continuous function so that with $r(t,\lambda)>0$ we have
\begin{align}\label{car_angle}
(1,-i)\cdot v(t,\lambda)=r(t,\lambda) e^{i \phi_\lambda(t)/2}, \qquad \mbox{ with }\phi_\lambda(0)=\varphi_0.
\end{align}
More precisely, $t\mapsto i\varphi_\lambda(t)/2+\log r(t,\lambda)$ is the unique lifting of the 2-dimensional curve $(1,-i)\cdot v$ under the covering given by the exponential function from $\CC$ to $\CC\setminus\{0\}$, with initial condition $i \varphi_0/2$. 
We will call $\phi_\lambda(t)$ the \emph{phase function} of $\tau$.

Some authors define the  phase function as $2 \phi_\lambda$. With our definition the phase angle has period $2\pi$ which is more convenient for us. In particular, $v(t,\lambda)\parallel u_1$ if and only if $\phi_\lambda(t)=\varphi_1$ mod $2\pi$.

If $v=(v_1,v_2)^t$ solves the ODE (\ref{car_angle}) then $r_\lambda(t)=\frac{v_1(t,\lambda)}{v_2(t,\lambda)}=-\cot(\phi_\lambda(t)/2)$ satisfies
\begin{align*}
r_\lambda'&=\frac{v_1' v_2-v_2' v_1}{v_2^2}=\frac{v^t J v'}{v_2^2}=\lambda \frac{v^t R v}{v_2^2}=\lambda (z,1) R \bin{z}{1}.
\end{align*}
If we now assume that $R$ is of the form (\ref{Rrep}) then we can further simplify this as
\begin{align}\label{carouselODE}
r_\lambda'=\lambda f \frac{y^2+(x-r_\lambda)^2}{y}, \qquad r_\lambda(0)=-\cot(\varphi_0/2).
\end{align}
For $\lambda>0$ the function $r_\lambda(t)$ is strictly increasing in $t$ and it restarts at $-\infty$ whenever it blows up to $\infty$. For $\lambda<0$ the function $r_\lambda(t)$ is strictly decreasing and restarts at $\infty$ after exploding to $-\infty$.

The ODE (\ref{carouselODE}) has a nice geometric representation. By (\ref{hyprot}) it describes the evolution of a boundary point $r_\lambda(t)\in \partial \HH$ which is continuously rotated with rate $2\lambda f(t)$ about the moving center of rotation $x(t)+i y(t)\in  \HH$.

The evolution of the angle $\phi_\lambda$ can be expressed using (\ref{carouselODE}) and $r_\lambda=-\cot(\phi_\lambda/2)$.
 One gets
\begin{align}\label{car_disk}
\phi_\lambda'=2f \lambda \frac{\left|e^{i \phi_\lambda}- \gamma\right|^2}{1-| \gamma|^2}, \qquad \phi_\lambda(0)=\varphi_0.
\end{align}
Here $ \gamma=U(x+i y)\in \UU$ is the representation of the path $x+i y$ in the Poincar\'e disk model (see (\ref{Cayley})). One advantage of this representation is that the solution has no blow-ups. The geometric picture is the same as before, but now in the Poincar\'e disk: the boundary point $e^{i \phi_\lambda(t)}\in \partial \UU$ is   continuously rotated with rate $2\lambda f(t)$ about the moving center of rotation $\gamma(t) \in  \UU$. Note that Proposition \ref{prop:ODE} implies that $\phi_\lambda(t)$ is analytic in $\lambda$ for any fixed $t\in [0,T)$.

\subsection{End behavior of the phase function}

The following theorem fully describes the behavior of the phase function as $t\to T$ and shows how it can be used to describe the eigenvalues in the limit circle case.

\begin{theorem}\label{thm:osc}
Consider a Dirac operator $\tau=\Dirop(x+i y, u_0,u_1, f)$ satisfying the conditions (\ref{cond1})-(\ref{cond3}), with non-parallel boundary conditions $u_i\parallel (\cos(\varphi_i/2), -\sin (\varphi_i/2))^t$, and a Hilbert-Schmidt inverse.
Let $\phi_\lambda(t)$ be the phase function introduced in (\ref{car_angle}). Then for every $\lambda$,
\[
\phi_\lambda(T)=\lim_{t\uparrow  T} \phi_\lambda(t)
\]
exists and the limit is finite.

In the limit circle case the function $\lambda\to \phi_\lambda(T)$ is continuous and strictly increasing, and
\begin{align*}
\{\textup{eigenvalues of } \tau\}=\{\lambda: \phi_\lambda(T)\in  \varphi_1+2\pi \mathbb Z\}.
\end{align*}
In the limit point case we have $\phi_\lambda(T)\in \varphi_1+2\pi \Z$ for all nonzero $\lambda\in\R$.
\end{theorem}

The theorem shows that in the limit circle case the endpoint of the carousel is a continuous, strictly increasing function of $\lambda$, while in the limit point case the endpoint (as a point on the unit circle) is always  $e^{i \varphi_1}$ for all $\lambda\not=0$.

\begin{proof}[Proof of Theorem \ref{thm:osc}]
The phase angle $\phi_\lambda(t)$ satisfies the ODE (\ref{car_disk})
where $\gamma$ is the image of the path $x+i y$ under the transformation $U(z)=\frac{z-i}{i+z}$. Note that we have
\begin{align}\label{XXXb}
|\gamma|< 1, \qquad \frac{1}{1-|\gamma|^2}=\frac{x^2+(y+1)^2}{4 y}.
\end{align}
This implies
\begin{align}\label{XXXc}
c_0 \|R\|\le \frac{f}{1-|\gamma|^2}\le c_1 \|R\|,
\end{align}
i.e.~the function  $ \frac{f}{1-|\gamma|^2}$ is integrable on $[0,T)$ if and only if $\|R\|$ is integrable there.

Note that for $\lambda=0$ we have $\phi_0(t)=\varphi_0$ for all $t$. We will now assume that $\lambda>0$, the $\lambda<0$ case can be treated similarly. For $\lambda>0$ the function  $\phi_\lambda(t)$ is increasing in $t$, which shows  that the limit $\lim_{t\uparrow T}\phi_\lambda(t) $ exists. To show that the limit is finite we need to prove that $\phi_\lambda(t)$ is bounded on $[0,T)$ for each $\lambda$.

We first consider the limit circle case. As $\left|e^{i \phi_\lambda}- \gamma\right|\le 2$, we have
\begin{align*}
\phi_\lambda(T)-\phi_\lambda(t) \le  8\lambda \int_t^T \frac{f}{1-|\gamma|^2}ds.
\end{align*}
Since $\|R\|$ is integrable on $[0,T)$ this means that $\frac{f}{1-|\gamma|^2}$ is also integrable. This shows that the limit is finite, the convergence is uniform on compact sets of $\lambda$, and hence the limit $\lambda\to \phi_\lambda(T)$ is continuous. The geometric picture behind ODE (\ref{car_disk}) shows that $\phi_\lambda(t)$ is strictly increasing as a function of $\lambda$ for any $t>0$. This also implies that $\phi_\lambda(T)$ is non-decreasing in $\lambda$. The fact that it is strictly increasing follows from the fact that if $\phi_{\lambda_1}(T)=\phi_{\lambda_2}(T)$ for some $\lambda_1<\lambda_2$ then the ODE (\ref{car_disk}) would imply that $\phi_{\lambda_1}(t)>\phi_{\lambda_2}(t)$ for a $t$ close to $T$ which contradicts the fact that $\lambda\to \phi_\lambda(t)$ is strictly increasing.

The integrability of $\|R\|$  also implies that the solution of (\ref{osc}) can be extended to $[0,T]$ in a continuous way, and it  is in $L^2_R$. The eigenvalues of $\tau$ are exactly the $\lambda\in \R$ values for which the shooting problem $\tau v=\lambda v$, $v(0)=u_0$, $v(T)\parallel u_1$ can be solved. But this is equivalent to $
\phi_\lambda(T)\in  \varphi_1+2\pi \mathbb Z$.

Getting back to the limit point case, let us consider  an approximating sequence $T_n\uparrow T$ with the corresponding operators $\tau_n$. Then  the eigenvalues of $\tau_n$ converge to those of $\tau$ as $n\to \infty$. Since $\tau_n$ is limit circle, the number of eigenvalues of $\tau_n$ in $[0,\lambda]$ is given by $\left|\{\lambda: \phi_\lambda(T_n)\in  \varphi_1+2\pi \mathbb Z\}\right|$. This gives an upper bound on $\limsup_{n\to \infty} \phi_\lambda(T_n)$ in terms of the number of eigenvalues of $\tau$ in $[0,\lambda]$, which shows that $\lim_{t\uparrow T}\phi_\lambda(t) $ is finite.

%
%

In the limit point case we have $\int_0^T \|R\| dt=\infty$. By \eqref{XXXc} we have $\int_0^T \frac{f}{1-|\gamma|^2} dt=\infty$.
%
%
Using   equation (\ref{car_disk}) and the finiteness of $\phi_\lambda(T)$  we see that
$
f \frac{\left|e^{i \phi_\lambda}-\gamma\right|^2}{1-|\gamma|^2}
$
is integrable on $[0,T)$. Note, that for the vector $a=(e^{i\alpha/2},e^{-i\alpha/2})^t$, and $\tilde U=\tfrac{1}{\sqrt{2}} \mat{1}{-i}{1}{i}$ we have
$$
(\tilde U^*a)^*R\tilde U^*a=\frac{ f |e^{i\alpha}-\gamma|^2}{1-|\gamma|^2}.
$$
By definition of the limit point case this is integrable if and only if $\tilde U^*a$ is parallel to $u_1$, which is equivalent to  $\alpha\in \varphi_1+2\pi \Z$. Let $\alpha$ be such an angle, and let $d=|e^{i\alpha}-e^{i\phi_\lambda(T)}|/2$. By the triangle inequality
$$|e^{i\alpha}-\gamma(t)|+|\gamma(t)-e^{i\phi_\lambda(t)}|+|e^{i\phi_\lambda(t)}-e^{i\phi_\lambda(T)}|\ge 2d.$$
The inequality $(x+y+z)^2\le 3(x^2+y^2+z^2)$ gives
$$|e^{i\alpha}-\gamma(t)|^2+|\gamma(t)-e^{i\phi_\lambda(t)}|^2+|e^{i\phi_\lambda(t)}-e^{i\phi_\lambda(T)}|^2\ge \frac{4}{3}d^2.$$
We consider two cases depending on whether the last term is at most $d^2$ or more. If it is more, 4 is still an upper bound. This gives
\begin{align}\label{xyzd}
|e^{i\alpha}-\gamma|^2+|\gamma-e^{i\phi_\lambda(t)}|^2+4\cdot \mathbf{1}(|e^{i\phi_\lambda(t)}-e^{i\phi_\lambda(T)}|>d)\ge \tfrac{d^2}{3}.
\end{align}
Multiplying  inequality (\ref{xyzd}) with $\frac{f}{1-|\gamma|^2}$ and integrating on $[0,T)$, we see that all three
integrals on the left hand side are finite. But the right side can only be finite if $d=0$,   showing that  $\phi_\lambda(T)\in \varphi_1+2\pi \Z$.
\end{proof}

\subsection{The hyperbolic carousel}

The evolution (\ref{car_disk})  has already appeared in \cite{BVBV} where it was used to define the hyperbolic carousel. This is a  geometric functional producing a discrete set of points from a hyperbolic path and two boundary points.

\begin{definition}[Hyperbolic carousel]\label{def:hypc}
Suppose that $\gamma(t), t\in [0,T)$ is a measurable, locally bounded path in the hyperbolic plane $\HH$,   $\eta_0, \eta_1\in \partial \HH$ are two distinct boundary points and  $f$ is a positive, locally integrable  function on $[0,T)$. The hyperbolic carousel
associated to  $(\gamma, \eta_0, \eta_1, f)$ produces a discrete set of points on $\R$, denoted by $\mathcal{HC}(\gamma, \eta_0, \eta_1, f)$, defined as follows.

For any fixed $\lambda\in \R$ we  consider a moving boundary point $r_\lambda(t)$ with $r_\lambda(0)=\eta_0$, which is rotated about $\gamma(t)$ continuously with rate $2 \lambda f(t)$ for $t\in [0,T)$.
 For each $\lambda$ we count how many times the moving point passes $\eta_1$ (counting it with a negative sign if $\lambda<0$) and we denote this number by $N(\lambda)$.
 More precisely, $N(\cdot)$ is the right-continuous version of the function
 $$\lambda \mapsto\operatorname{sgn}(\lambda) \cdot \left|\left\{t\in[0,T): r_\lambda(t)=\eta_1\right\}\right|.$$
  If $N(\cdot)$ is a finite function then $\mathcal{HC}(\gamma, \eta_0, \eta_1, f)$ is the set of points whose counting function is $N$. If  $N(\cdot)$ is not finite then $\mathcal{HC}(\gamma, \eta_0, \eta_1, f)$ is undefined.
  \end{definition}

 If $f$ is the constant function $\tfrac12$ then we drop it from the notation, i.e.~$\mathcal{HC}(\gamma, \eta_0, \eta_1)=\mathcal{HC}(\gamma, \eta_0, \eta_1, \tfrac12)$. This corresponds to the carousel where the boundary point $r_\lambda$ is rotated with constant speed $\lambda$.

Note that the definition does not rely on any particular representation of the hyperbolic plane.
If we consider the Poincar\'e disk representation of $\HH$ then we can describe the moving boundary points of the carousel as $r_\lambda(t)=e^{i \phi_\lambda(t)}$ with a continuous $ \phi_\lambda(t)$ satisfying $\phi_\lambda(0)=\arg \eta_0$. By the discussion around (\ref{car_disk}) we see that the function $ \phi_\lambda(t)$ satisfies the ODE (\ref{car_disk}) with initial condition  $\varphi_0=\arg \eta_0$. Thus the ODE describing the moving boundary points of the carousel is exactly the same as the ODE for the phase angle of a Dirac operator. The next proposition shows that under some mild conditions the set of points produced by the carousel is exactly the same  as the spectrum of the Dirac operator built from the same ingredients.

\begin{proposition}\label{prop:car_dir}
Assume that $u_0, u_1$ are nonzero vectors that are not parallel, the operator $\Dirop(x+i y, f, u_0, u_1)$ defined on $[0,T)$ satisfies conditions (\ref{cond1})-(\ref{cond3}) and this operator has a Hilbert-Schmidt inverse. Then $\mathcal{HC}(x+i y, u_0, u_1, f)$ is well defined and it is equal to the spectrum of $\Dirop(x+i y,  u_0, u_1,f)$.
\end{proposition}

\begin{proof}
Let $\varphi_0, \varphi_1\in [0,2\pi)$ so that $u_i \parallel (\cos(\varphi_i/2), -\sin(\varphi_i/2))^t$. We have seen that the phase angle of  $\tau=\Dirop(x+i y, u_0, u_1,f)$ satisfies the ODE (\ref{car_disk}). This phase angle also encodes the moving boundary points of the carousel with the parameters $(x+iy, u_0, u_1, f)$. This means that the counting function of $\mathcal{HC}(x+i y, u_0, u_1, f)$ is given by the right-continuous version of the function
\begin{align*}
\tilde N(\Lambda) &=\operatorname{sign}(\Lambda) \cdot \left| \{t\in [0,T): \phi_\Lambda(t)\in \varphi_1+2\pi\Z \} \right|.
\end{align*}
In order to prove the proposition we need to show that this is finite for all $\Lambda$ and that the right-continuous version of $\tilde N(\cdot)$ is exactly the counting function of the spectrum of $\tau$.

For simplicity, we will only deal with the $\Lambda>0$ case (the other case can be handled the same way).

Consider an approximating sequence $T_n\uparrow T$ with the corresponding operators $\tau_n$, as in Theorem \ref{thm:appr}.
Denote the counting function of the spectrum of $\tau$ and $\tau_n$ by $F$ and $F_n$, respectively.
Fix $n$.  By Theorem \ref{thm:osc} we have
\begin{align*}
F_n(\Lambda)&= \left|\{\lambda\in [0,\Lambda]: \phi_\lambda(T_n)\in \varphi_1  +2\pi\Z\}\right|.
\end{align*}
From the geometric definition (or the ODE (\ref{car_disk})) it follows that for a fixed $\lambda>0$ the function $t\to \phi_\lambda(t)$ is strictly increasing, and for any fixed $0<t<T$ the function $\lambda\to \phi_\lambda(t)$ is also strictly increasing. Note that $\phi_\lambda(t)$ is continuous for $(\lambda,t)\in [0,\Lambda]\times[0,T)$, and we have
  \[
  \phi_\lambda(0)=\phi_0(t)=\varphi_0, \qquad \textup{and} \qquad \varphi_0 \not\equiv \varphi_1 \mod 2\pi.
  \]
    For any $0<T_n<T$ the functions  $\lambda\to \phi_\lambda(T_n), \lambda\in [0,\Lambda]$ and $t\to \phi_\Lambda(t), t\in [0,T_n]$ are both continuous, strictly increasing and have the same starting and end points which yields
\begin{align}\label{car_123}
\left|\{\lambda\in [0,\Lambda]: \phi_\lambda(T_n)\in \varphi_1 + 2\pi\Z\}\right|=\left| \{t\in [0,T_n): \phi_\Lambda(t)\in \varphi_1 + 2\pi\Z\}\right|.
\end{align}
For any fixed $\Lambda$ the right side of (\ref{car_123}) converges to $\tilde N(\Lambda)$ as $n\to \infty$. The left side of  (\ref{car_123}) is $F_n(\Lambda)$ and by Theorem \ref{thm:appr} this will converge to $F(\Lambda)$ as $n\to \infty$ for every $\Lambda$ which is not an eigenvalue of $\tau$. From this it follows that $\tilde N(\cdot)$ is a finite function and that its right-continuous version is exactly $F(\cdot)$, the counting function of the spectrum of $\tau$.
\end{proof}

We have seen the carousel ODE in both the half-plane (\ref{carouselODE}) and the unit disk (\ref{car_disk}) coordinates. We have also seen that the driving path in the half-plane representation can be used to express the parameters of the corresponding Dirac operator. With a simple linear transformation of the operator  we can also recover the driving path in the unit disk coordinates. Suppose that $\tau=R(t) J \partial_t$ where $R(t)$ is given as in (\ref{Rrep}). Now consider the operator
\begin{align}\label{HtoU}
\tilde \tau u(t)=\tilde U \tau(\tilde U^{-1} u)=\tilde U R^{-1}(t) J \tilde U^{-1} u'(t)
\end{align}
 defined on functions $u:[0,1]\to \CC^2$ where $\tilde U=\tfrac{1}{\sqrt{2}} \mat{1}{-i}{1}{i}$. (This is just the linear transformation corresponding to the Cayley map $U(z)$.) A simple computation shows that
\begin{align}\label{Dirac_disk}
\tilde \tau u=\frac{1}{f(t)(1-|\gamma(t)|^2)} \mat{1+|\gamma(t)|^2}{2 \gamma(t)}{2 \bar \gamma(t)}{1+|\gamma(t)|^2} \mat{-i}{0}{0}{i} u'
\end{align}
where $\gamma(t)=U(x(t)+i y(t))$ is exactly the image of the driving path in the unit disk representation.

\subsection{The reverse phase function}

The phase function can be started at the right end point $T$, even when the operator is limit point there.

\begin{lemma}
Consider a Dirac operator $\tau=\Dirop(x+iy, u_0, u_1, f)$ on $[0,T)$ satisfying conditions (\ref{cond1})-(\ref{cond3}) with non-parallel boundary conditions and with a Hilbert-Schmidt inverse.
Let $\varphi_i\in [0,2\pi)$ be the angles with $u_i \parallel (\cos(\varphi_i/2), -\sin (\varphi_i/2))^t$ and let $\gamma=U(x+i y)$ be the representation of the path $x+i y$ in the Poincar\'e disk, as in (\ref{car_disk}).

Then there is a unique solution of the following ODE system
\begin{align}\label{phase_reversed}
\rho_\lambda'=2f \lambda \frac{\left|e^{i \rho_\lambda}- \gamma\right|^2}{1-| \gamma|^2}, \qquad \lim_{t\uparrow T}\rho_\lambda(t)=\varphi_1
\end{align}
with the following condition:  for $\lambda>0$ if the function $\tilde \rho_\lambda$ also solves (\ref{phase_reversed}) then $\rho_\lambda(t)\ge \tilde \rho_\lambda(t)$ for $t\in [0,T)$. Similarly: if $\lambda<0$ then  $\rho_\lambda(t)\le \tilde \rho_\lambda(t)$ for $t\in [0,T)$.

Moreover, $\lambda\to \rho_\lambda(t)$ is continuous and strictly decreasing for any $t\in [0,T)$, and
\begin{align}\label{rev_spec}
\{\textup{eigenvalues of } \tau\}=\{\lambda: \rho_\lambda(0)\in  \varphi_0+2\pi \mathbb Z\}.
\end{align}
\end{lemma}

We  call the function $\rho_\lambda$ the \emph{reverse phase function} of $\Dirop(x+iy, u_0, u_1, f)$.
\begin{proof}
If the function $\|R\|$ is bounded on $[0,T)$ then the lemma follows  from the
 time reversal $t\to T-t$ and Theorem \ref{thm:osc}.

Set $0<\tilde T<T$. Since $\|R\|$ is bounded on $[0,\tilde T]$, the reverse phase function $ \rho_{\lambda, \tilde T}$ for the restriction of $\tau$ to $[0,\tilde T)$ exists, and satisfies  the ODE
\begin{align*}
 \rho_{\lambda, \tilde T}'=2f \lambda \frac{\left|e^{i \rho_{\lambda, \tilde T}}- \gamma\right|^2}{1-| \gamma|^2}, \qquad \rho_{\lambda, \tilde T}(\tilde T)=\varphi_1.
\end{align*}
The function  $\lambda\to \rho_{\lambda,\tilde T}(t)$ is continuous and strictly increasing for any $t\in [0,\tilde T)$, and the spectrum of $\tau$ restricted to $[0,\tilde T)$ is equal to the set  $\{\lambda: \rho_{\lambda, \tilde T}(0)\in  \varphi_0+2\pi \mathbb Z\}$.

We will show that the limit $\lim_{\tilde T\uparrow T} \rho_{\lambda,\tilde T}$   satisfies the ODE (\ref{phase_reversed}) and  the listed conditions. Without loss of generality we assume $\lambda>0$.

The ODE (\ref{phase_reversed})
describes the evolution of $e^{i \rho_{\lambda,\tilde T}}$ as a boundary point of $\partial \UU$  being continuously rotated about $\gamma$ with rate $2f\lambda$. This implies the following coupling result: if $0<T_1<T_2<T$ then $\rho_{\lambda, T_2}(t)< \rho_{\lambda, T_1}(t)$ for  $t\in [0,T_1]$. We also have $\rho_{0,T_1}(t)=\rho_{0,T_2}(t)=\varphi_1$. Borrowing the arguments of the proof of Theorem \ref{thm:osc} we see that the limit
\begin{align}\label{rho}
\rho_\lambda(t):= \lim_{\tilde T\uparrow T} \rho_{\lambda,\tilde T}(t)
\end{align}
 exists for each $t\in [0,T)$ because of  monotonicity, and $\rho_\lambda(t)$ will be non-decreasing in $t$.
  For a fixed $\Lambda>0$ the limit $|\rho_{\lambda}(0)|$ can be bounded uniformly on $[0,\Lambda]$  in terms of the number of eigenvalues of $\tau$ in $[0,\Lambda]$. This shows that the limit $\rho_\lambda(0)$ is finite, which in turn shows that $\rho_\lambda(t)$ is finite for any $t\in [0,T)$. Note that since $\lambda\to \rho_{\lambda,\tilde T}(t)$ is decreasing, the function $\lambda\to \rho_{\lambda}(t)$ is non-increasing.

For any $0\le t_0<t_1\le \tilde T$ we have
\[
 \rho_{\lambda,\tilde T}(t_1)- \rho_{\lambda,\tilde T}(t_0)=\int_{t_0}^{t_1} 2f \lambda \frac{\left|e^{i \rho_{\lambda, \tilde T}(u)}- \gamma(u)\right|^2}{1-| \gamma(u)|^2}du.
\]
Since $\gamma$ is locally bounded in $\UU$ on $[0,T)$ and $f$ is locally bounded on $[0,T)$ the integral on the right converges to $\int_{t_0}^{t_1} 2f \lambda \frac{\left|e^{i \rho_{\lambda}(u)}- \gamma(u)\right|^2}{1-| \gamma(u)|^2}du$, which shows that the limit $ \rho_{\lambda}$  satisfies the ODE (\ref{phase_reversed}) on $[0,T)$.

We will now prove the continuity of the function $\lambda\to \rho_\lambda(0)$. Once we have that, the continuity of $\lambda\to \rho_\lambda(t)$ for $t\in (0,T)$ follows by the continuous dependence of the solution of an ODE on the initial parameters.

The limit circle case can be handled the same way as it was done in the proof of Theorem \ref{thm:osc}. There is a unique solution of the ODE (\ref{phase_reversed}) which can be extended continuously to $[0,T]$ with $\rho_\lambda(T)=\varphi_1$. The solution $\rho_\lambda(t)$ is continuous on $(\lambda,t)\in \R\times [0,T]$, and it is strictly decreasing in $\lambda$ for any fixed $t$. Since $\rho_{\lambda,\tilde T}(0)$ is continuous, strictly increasing and converges to $\rho_{\lambda}(0)$ with these same properties we immediately get that the set $\{\lambda: \rho_{\lambda, \tilde T}(0)\in  \varphi_0+2\pi \mathbb Z\}$ converges pointwise to the set $\{\lambda: \rho_{\lambda}(0)\in  \varphi_0+2\pi \mathbb Z\}$,  proving (\ref{rev_spec}).

Consider now the limit point case and assume that there is a $\lambda_0$ where $ \rho_\lambda(0)$ is not continuous. Since $ \rho_\lambda(0)$ is non-increasing in $\lambda$ this means that the left and right limits $\rho_{\lambda_0^-}(0), \rho_{\lambda_0^+}(0)$ at $\lambda_0$  exist and $\rho_{\lambda_0^+}(0)<\rho_{\lambda_0^-}(0)$. Choose  $\xi\neq \varphi_1$ so that $\rho_{\lambda_0^+}(0)<\xi <\rho_{\lambda_0^-}(0)$. Then using the definition of $\rho_\lambda$ as a limit we see that for small enough $\eps>0$ there is a  $\tilde T_0>0$ so that for $\tilde T>\tilde T_0$ we have
\[
\rho_{\lambda_0+\eps, \tilde T}(0)<\xi <\rho_{\lambda_0-\eps, \tilde T}(0)
\]
This means that if we consider the operator $\tau$, but with initial condition $\tilde u_0=(\cos(\xi/2),- \sin(\xi/2))^t$, then the approximating operators on $[0,\tilde T)$ will all have an eigenvalue in $[\lambda-\eps, \lambda+\eps]$. By  Theorem \ref{thm:appr} this will also be true for $\tau$ itself. Since this is true for any small enough $\eps$, this means that $\lambda$ is an eigenvalue for $\tau$ with initial condition  $\tilde u_0$. But $\xi$ is arbitrary from $(\rho_{\lambda_0^+}(0),\rho_{\lambda_0^-}(0))$, thus there are at least two linearly independent solutions of $\tau v=\lambda v$ that are in $L^2_R$ near the right boundary $T$.  By Theorem \ref{thm:Weyl} this would imply that $\tau$ is actually limit circle at $T$, and this contradiction proves the continuity of $\rho_\lambda(0)$.

Next we show that $\rho_\lambda(0)$ is strictly decreasing in $\lambda$. Assume that $\rho_{\lambda_1}(0)=\rho_{\lambda_2}(0)$ for some $\lambda_1<\lambda_2$. Then from the ODE (\ref{phase_reversed})  it follows that there is an $\eps>0$ so that $\rho_{\lambda_1}(\eps)<\rho_{\lambda_2}(\eps)$. But this contradicts the fact that $\lambda\to \rho_{\lambda}(\eps)$ is  non-increasing, which we have seen already.

This shows that $\lambda\to \rho_\lambda(0)$ is strictly decreasing and continuous. This is also true for $\lambda\to \rho_{\lambda, \tilde T}(0)$ for each $\tilde T$, so by (\ref{rho}) we have that the sets $\{\lambda: \rho_{\lambda, \tilde T}(0)\in  \varphi_0+2\pi \mathbb Z\}$ converge to the set $\Xi=\{\lambda: \rho_{\lambda}(0)\in  \varphi_0+2\pi \mathbb Z\}$. But by Theorem \ref{thm:appr} this implies that the spectrum of $\tau$ is given  by the set $\Xi$, proving (\ref{rev_spec}).

Finally, we have to show that $\rho_\lambda$ is the unique solution of (\ref{phase_reversed}) satisfying the second (minimal/maximal) condition. It is enough to show that $\rho_\lambda$ has the prescribed property, the uniqueness follows (since we can only have one minimal or maximal solution.) Suppose  $\tilde \rho_\lambda$ is another solution of (\ref{phase_reversed}) with $\rho_\lambda(t)<\tilde \rho_\lambda(t)$ for some $t$. Using the definition (\ref{rho}) we can find $t<\tilde T$ so that $\rho_{\lambda, \tilde T}(t)<\tilde \rho_\lambda(t)$. For a fixed $\lambda$ the solutions of our ODE do not cross, which gives  $\varphi_1=\rho_{\lambda, \tilde T}(\tilde T)<\tilde \rho_\lambda(\tilde T)$.
Since $\tilde \rho_\lambda(u)$ is increasing in $u$, we have $\tilde \rho_\lambda(\tilde T)<\lim_{u\to T}\tilde \rho_\lambda(u)=\varphi_1$, which is a contradiction.
\end{proof}

\section{Unitary matrices as  Dirac operators}\label{s:unitary}

The goal of this section is to show that a finite unitary matrix can be connected to a hyperbolic carousel with a piecewise constant driving path, and consequently to a Dirac operator with a piecewise constant weight  function $R(t)$.

\subsection{The Szeg\H o recursion} \label{subs:Szego}

Consider an  $n\times n$ unitary matrix $V$ with $n$ distinct eigenvalues. Fix a unit vector $e$ which is not orthogonal to any of the eigenvectors, then the
vectors $e, Ve, \dots, V^{n-1} e$
form a basis.  Applying the Gram-Schmidt procedure yields an orthogonal basis
\begin{align}
\Phi_0(V) e, \Phi_1(V) e, \dots, \Phi_{n-1}(V) e, \label{GSbasis}
\end{align}
where $\Phi_0(z)=1$ and $\Phi_k(z)$ is a monic polynomial of degree $k$ for $k>0$. This sequence can be naturally extended for $k=n$, with $\Phi_n(z)$ defined as $\det(z-V)$, the characteristic polynomial of $V$.
Together with the reversed polynomials  $\Phi_k^*(z)=z^k \overline{\Phi_k(1/\bar z)}$, they satisfy the famous Szeg\H o recursion (see e.g.~Section 1.5 of \cite{OPUC1}):
\begin{align}\label{Szego1}
\bin{\Phi_{k+1}}{\Phi_{k+1}^*}=A_k Z \bin{\Phi_{k}}{\Phi_{k}^*}, \qquad \bin{\Phi_0}{\Phi_0^*}=\bin{1}{1},\qquad 0\le k\le n-1,
\end{align}
where
$$
A_k=\mat{1}{-\bar \alpha_k}{-\alpha_k}{1}, \qquad  Z=\mat{z}{0}{0}{1}.
$$
The complex numbers $\alpha_0, \dots, \alpha_{n-1}$ are called Verblunsky coefficients.  They satisfy $|\alpha_k|<1$ for  $0\le k<n-1$, and $|\alpha_{n-1}|=1$. The Verblunsky coefficients determine the recursion which in turn can be used to identify the eigenvalues of $V$ as the roots of $\Phi_n(z)$.

By expanding the last step in the recursion, we see that  $\Phi_n(z)=0$ if and only if
\begin{align}\label{end_cond}
Z\bin{\Phi_{n-1}(z)}{\Phi_{n-1}^*(z)} \parallel \bin{\bar \alpha_{n-1}}{1}.  
\end{align}
Although the Verblunsky coefficients cannot identify the matrix $V$, they determine its spectral measure corresponding to the vector $e$ and vice versa. In \cite{CMV} the authors construct a 5-diagonal unitary matrix in terms of the Verblunsky coefficients which is similar to $V$. The constructed matrix is called the CMV representation of $V$, see \cite{Simon_CMV} for additional details.

\subsection{Operator from the Szeg\H o recursion }

The recursion (\ref{Szego1}) for $(\Phi_k(z), \Phi_k^*(z))^t$ depends on $z$ via the matrix $Z$, the part involving the $A_k$ matrices is the same for each $z$.  The goal of this section is to separate these two components, and to show that this leads to a Dirac operator and a hyperbolic carousel.

Write  $z=e^{i \lambda}$ for $\lambda\in \mathbb C$ (note that $\lambda \in \R$ is the  most relevant case). Set
\begin{align*}
f_0=\bin{1}{1}, \quad\textup{and} \quad f_{k+1}=f_{k+1}(\lambda)=e^{-i
\lambda (k+1)/2} M_{k}^{-1} Z\bin{\Phi_{k}(e^{i \lambda})}{\Phi_{k}^*(e^{i \lambda})}, \qquad 0\le k\le n-1,
\end{align*}
where $M_k=A_{k-1}\cdots A_0$ with $M_0=I$.
Then  the sequence $f_k$ satisfies the recursion
\begin{align}\label{frec}
f_{k+1}=\mat{e^{i\lambda/2}}{0}{0}{e^{-i\lambda/2}}^{M_{k}} f_k, \qquad f_0=\bin{1}{1},\qquad 0\le k\le n-1,
\end{align}
where  $X^Y=Y^{-1}X Y$ is our notation for conjugation. By  \eqref{end_cond}, $z=e^{i \lambda}$ is an eigenvalue if and only if \[ f_n(\lambda) \parallel M_{n-1}^{-1} \bin{\bar \alpha_{n-1}}{1}.\]

Now let $M(t)=M_{\lfloor nt \rfloor}, t\in [0,1)$. Consider the differential operator $\tau$ acting on functions $g: [0,1)\to \CC^2$ as
\begin{align}\label{tauXXX}
\tau g=2 \mat{-i}{0}{0}{i}^{M(t)} g',
\end{align}
with initial and end conditions
\begin{align}\label{end_cond_x}
g(0)\parallel u_0:=\bin{1}{1}, \qquad g(1)\parallel u_1:=M_{n-1}^{-1} \bin{\bar \alpha_{n-1}}{1}.
\end{align}
The solution of the eigenvalue equation $\tau g=\mu g$ satisfies
\[
g'(t)= \mat{i\mu/2}{0}{0}{-i\mu/2}^{M(t)} g(t).
\]
As $M(t)$ is constant on intervals of the form $[\frac{k}{n},\frac{k+1}{n})$, we can explicitly solve the ODE to get
\[
g(\tfrac{k+1}{n})=\mat{e^{i\tfrac{\mu}{2n}}}{0}{0}{e^{-i\tfrac{\mu}{2n}}}^{M(k/n)} g(\tfrac{k}{n}).
\]
Recalling $M(k/n)=M_k$ and the recursion (\ref{frec}) we get that $g(k/n)=g_\mu(k/n)=f_k(\tfrac{\mu}{n})$, and that the eigenvalues of $\tau$ are given by the set
\begin{align}\label{yyy_discrete}
\{\mu\in \mathbb R:  e^{i\mu/n} \textup{ is an eigenvalue of the Szeg\H o recursion}\}.
\end{align}
Note that $M_k=\mat{p_k}{q_k}{\bar q_k}{\bar p_k}$ with $|p_k|^2-|q_k|^2>0$,  since the $A_k$ matrices are of this form, and this property is inherited in products. A simple computation shows that for $0\le k\le n-1$ we have
\begin{align}\label{bbrec}
\mat{-i}{0}{0}{i}^{M_k}=\frac{1}{1-|b_k|^2} \mat{1+|b_k|^2}{2 b_k}{2 \bar b_k}{1+|b_k|^2} \mat{-i}{0}{0}{i}, \quad b_k=-\frac{q_k}{p_k}=\mathcal P M_k^{-1} \bin{0}{1}.
\end{align}
This means that the differential operator $\tau$ defined in (\ref{tauXXX}) is exactly of the form of (\ref{Dirac_disk}) with $f=\frac12$. Thus it is just a linear conjugate of a Dirac operator of the form (\ref{Diracop}). The corresponding function $R_t$  will be piecewise constant, so the resulting Dirac operator is limit circle at 1.

Note that
\begin{align}\label{b_k}
b_k=\mathcal P A_0^{-1}\cdots A_{k-1}^{-1} \bin{0}{1}=\mathcal P \mat{1}{\bar \alpha_0}{\alpha_0}{1}\cdots \mat{1}{\bar \alpha_{k-1}}{\alpha_{k-1}}{1} \bin{0}{1},
\end{align}
for $k\le n-1$. Extending this definition to $k=n$ and comparing it to (\ref{end_cond_x}) we get that $|b_n|=1$ and $u_1\parallel (b_n, 1)^t$.

Thus we have shown the following.

\begin{proposition}\label{prop:unitary_car}
Suppose that $\alpha_0, \dots, \alpha_{n-1}$ are the Verblunsky coefficients of an $n\times n$ unitary matrix with distinct eigenvalues $\{e^{i\lambda_k}, 1\le k\le n\}$.
Let $b_0=0$ and define
\begin{align}\label{car_center}
b_k=\mathcal P \mat{1}{\bar \alpha_0}{\alpha_0}{1}\cdots \mat{1}{\bar \alpha_{k-1}}{\alpha_{k-1}}{1} \bin{0}{1},
\quad b(t)=b_{\lfloor nt \rfloor}, t\in [0,1].
\end{align}
Then $|b(t)|<1$ for $t\in [0,1)$ and $|b(1)|=1$. The Dirac operator
\begin{align}\label{unop}
\tau g=\frac{2}{1-|b|^2} \mat{1+|b|^2}{2 b}{2 \bar b}{1+|b|^2} \mat{-i}{0}{0}{i} g', \qquad t\in [0,1)
\end{align}
acting on functions $g: [0,1)\to \CC^2$ with initial and end  conditions $(1,1)^t$ and $(b(1),1)^t$ is self-adjoint on the appropriately defined domain. The spectrum of $\tau$ is
given by the set $\{n \lambda_k+2\pi n j: 1\le k\le n, j\in \ZZ\}$.
\end{proposition}
Equivalently, this set is obtained from the hyperbolic carousel $\mathcal{HC}(b, 1, b(1))$, where  the parameters are given the Poincar\'e disk coordinates.


\subsection{Operator from the modified Szeg\H o recursion}

The deformed (or modified) Verblunski coefficients were introduced in \cite{BNR2009} and are more natural in certain settings than the ordinary ones. Consider again an $n\times n$ unitary matrix with distinct eigenvalues and consider the orthogonal polynomials corresponding to a fixed  unit vector $e$ (not orthogonal to any of the eigenvectors), as discussed in Subsection \ref{subs:Szego}. Assume that 1 is not an eigenvalue, then the values $\Phi_k(1), \Phi_k^*(1)$ are all nonzero  and we can introduce the  polynomials
$$\Psi_k(z)=\frac{\Phi_k(z)}{\Phi_k(1)}, \qquad \Psi_k^*(z)=\frac{\Phi_k^*(z)}{\Phi^*_k(1)}.$$
A simple computation shows that $(\Psi_k, \Psi_k^*)^t$ satisfies exactly the same recursion \eqref{Szego1} as  $
(\Phi_k, \Phi_k^*)^t$, but with matrices
\begin{align}\label{tA_gamma}
\tilde A_k=\mat{\frac{1}{1-\gamma_k}}{-\frac{\gamma_k}{1-\gamma_k}}{-\frac{\bar \gamma_k}{1-\bar \gamma_k}}{\frac{1}{1-\bar \gamma_k}}, \qquad \gamma_k=\bar \alpha_k \frac{\Phi_k^*(1)}{\Phi_k(1)}
\end{align}
instead of $A_k$. The complex numbers   $\gamma_k$  are called deformed Verblunski coefficients. By \eqref{tA_gamma} we have $|\gamma_k|=|\alpha_k|$.  By Proposition 2.4 in \cite{BNR2009} the sequence of deformed Verblunski coefficients also determines the sequence of `regular' Verblunski coefficients.

From (\ref{end_cond}) and (\ref{tA_gamma})  we see that $z$ is an eigenvalue if and only if
\begin{align}\label{end_cond1}
Z\bin{\Psi_{n-1}(z)}{\Psi_{n-1}^*(z)} \parallel \bin{\gamma_{n-1}}{1}.
\end{align}
Repeating the arguments of the previous subsection we see that the following version of Proposition \ref{prop:unitary_car} holds.
\begin{proposition}\label{prop:unitary_car2}
Suppose that $\gamma_0, \dots, \gamma_{n-1}$ are the deformed Verblunsky coefficients of an $n\times n$ unitary matrix with eigenvalues $\{e^{i\lambda_k}, 1\le k\le n\}$. Define
\begin{align}\label{car_center2}
\tilde A_k=\mat{\frac{1}{1-\gamma_k}}{-\frac{\gamma_k}{1-\gamma_k}}{-\frac{\bar \gamma_k}{1-\bar \gamma_k}}{\frac{1}{1-\bar \gamma_k}}, \quad b_k=\mathcal P \tilde A_0^{-1}\cdots \tilde A_{k-1}^{-1}\bin{0}{1},
\quad b(t)=b_{\lfloor nt \rfloor}, \,\,\, t\in [0,1].
\end{align}
Then $|b(t)|<1$ for $t\in [0,1)$ and $|b(1)|=1$. The Dirac operator
\begin{align}\label{unop2}
\tau g=\frac{2}{1-|b|^2} \mat{1+|b|^2}{2 b}{2 \bar b}{1+|b|^2} \mat{-i}{0}{0}{i} g', \qquad t\in [0,1)
\end{align}
acting on functions $g: [0,1)\to \CC^2$ with initial and end conditions $(1,1)^t$ and $(b(1),1)^t$ is self-adjoint. The spectrum of $\tau$ is given by the set $\{n \lambda_k+2\pi n j: 1\le k\le n, j\in \ZZ\}$.
\end{proposition}
Note that the matrices $\tilde A_k$ are of the form $\mat{x}{1-x}{1-\bar x}{\bar x}$, in particular they have a common eigenvector $(1,1)^t$ with eigenvalue 1. This means that as M\"obius transformations they fix the point 1 on the boundary of the unit disk.  The product $\tilde A_{k-1}\cdots \tilde A_0$ will have the same property. As a result this product is just a simple function of $b_k$:
\[
\tilde A_{k-1}\cdots \tilde A_0=\mat{\frac{1}{1-b_k}}{-\frac{b_k}{1-b_k}}{-\frac{\bar b_k}{1-\bar b_k}}{\frac{1}{1-\bar b_k}}.
\]
Together with (\ref{car_center2}) this shows that $b_{k+1}$ has to be a function of $b_k$ and $\gamma_k$. Indeed, a quick computation gives the recursion
\begin{align}\label{b_rec}
b_{k+1} =\frac{b_k+\gamma_k\frac{1-b_k}{1-\bar b_k}}{1+\bar b_k \gamma_k \frac{
1-b_k}{1-\bar b_k}}, \qquad b_0=0.
\end{align}
The $\tilde A_k$ matrices (and their products) correspond to isometries of the hyperbolic plane in the disk model. The fact that $(1,1)^t$ is an eigenvector with eigenvalue 1 shows that these isometries fix the point $1\in \partial \UU$  on the boundary. If we consider the half-plane representation and move this fixed point to $\infty$ then the isometries become simple affine maps.
This observation suggests that the Dirac operator and the driving path will simplify in half-plane coordinates (when the point $1$ is sent to $\infty$). Recall the  Cayley map $U(z)$ from (\ref{Cayley}) mapping $\HH$ to $\UU$, mapping  $\infty$ to $1$, and the corresponding unitary matrix $\tilde U$.

Let $x_k+i y_k=U^{-1}(b_k)$ and define $x(t)+i y(t)=x_{\lfloor nt\rfloor}+i y_{\lfloor nt \rfloor}$. Then we can rewrite $\tau$ in the half-plane as
$$\sigma: f \mapsto  \tilde U^{-1}  \tau ( \tilde U f)$$
with the transformed boundary conditions
\[
f_0 \parallel \bin{1}{0}, \qquad f_1 \parallel \bin{x_n}{1}.
\]
By the discussion around equation (\ref{Dirac_disk}) we get that  $\sigma f = R(t)^{-1} J f'$ where $R(t)=\frac{1}{2y}\mat{1}{-x}{-x}{x^2+y^2}$.

Set
\begin{align}\label{AH}
W_k  =  \hat U^{-1} \tilde A_k  \hat U =\mat{1}{2\Im \tfrac{\gamma_k}{1-\gamma_k}}{0}{1+2\Re \frac{\gamma_k}{1-\gamma_k} }.
\end{align}
The matrix $W_k $  is in the affine group \eqref{affine1} of real matrices of the form
\begin{equation}\label{affine}
\mat{1}{-x}{0}{y}, \qquad y>0,
\end{equation}
and so are $X_k:=W_{k-1} \cdots W_0$. From (\ref{car_center2}) and \eqref{AH} we see that
\begin{equation}\notag x_k+i y_k = \mathcal P X_k^{-1}\bin{i}{1}.
\end{equation}
Note that for matrices $X$ of the form \eqref{affine} we have $\mathcal P X^{-1} (i,1)^t=x+i y$, so we can conclude that
$$X_k=\mat{1}{-x_k}{0}{y_k}.$$
The identity
 $ X_{k+1}=W_k X_k$ is short for the recursion
\begin{align}\label{xyrec}
x_{k+1}=x_k+v_k y_k, \qquad y_{k+1}=y_k(1+w_k), \qquad w_k=2\Re \frac{\gamma_k}{1-\gamma_k} , \quad v_k=-2 \Im \frac{\gamma_k}{1-\gamma_k},
\end{align}
with initial condition $x_0=0, y_0=1$.


The proposition below summarizes our findings.
\begin{proposition}
\label{prop:unitary_car2b}
Suppose that $\gamma_0, \dots, \gamma_{n-1}$ are the deformed Verblunsky coefficients of an $n\times n$ unitary matrix with distinct eigenvalues $\{e^{i\lambda_k}, 1\le k\le n\}$.
Suppose that $x_k, y_k$ solves the recursion
(\ref{xyrec}) with  $x_0=0, y_0=1$, and set $x(t)+i y(t)=x_{\lfloor nt \rfloor}+i y_{\lfloor nt \rfloor}$ for $t\in [0,1]$.
Then $y(t)>0$ for $t\in [0,1)$ and $y(1)=0$. Then the spectrum of the operator $\Dirop(x+iy, \infty, x(1))$ is given by the set $\{n \lambda_k+2\pi n j: 1\le k\le n, j\in \ZZ\}$.
\end{proposition}
In particular, this is the point process given by the hyperbolic carousel $\mathcal{HC}(x+i y, \infty , x(1))$.

\section{Circular ensembles} \label{s:circ}

The circular $\beta$-ensemble is a  random point process  with $n$ points on the unit circle. Using angles to describe the positions of the points the joint density is given by
\begin{align}\label{circdens}
\frac{1}{Z_{n,\beta}} \prod_{1\le j<k\le n} |e^{i \lambda_j}-e^{i \lambda_k}|^\beta, \qquad \lambda_k\in [-\pi,\pi).
\end{align}
Here $\beta>0$ and $Z_{n,\beta}$ is an explicitly computable positive constant. The joint density can be thought of as the Gibbs measure for $n$ charged particles confined to the unit circle and interacting via the two-dimensional Coulomb law. The parameter $\beta$ plays the role of inverse temperature.

For $\beta=2$ the distribution of the ensemble is the same as the distribution of the spectrum of a uniformly chosen $n\times n$ unitary matrix.  This distribution is called the circular unitary ensemble.
The connection to random matrices for general $\beta$ parameter was provided by \cite{KillipNenciu}.
\begin{theorem}[\cite{KillipNenciu}]\label{thm:KN}
Suppose that $\alpha_0, \dots, \alpha_{n-1}$ are independent, rotationally invariant random variables so that $|\alpha_k|^2$ has $\rm{Beta}(1,\tfrac{\beta}{2}(n-k-1))$ distribution for $0\le k\le n-2$ and $|\alpha_{n-1}|=1$. Then the $n\times n$ CMV matrix with Verblunsky coefficients  $\alpha_0, \dots, \alpha_{n-1}$ has joint eigenvalue density given by  (\ref{circdens}).
\end{theorem}
In \cite{KS} the authors used this representation to derive the point process limit of the circular $\beta$-ensemble, see Theorem \ref{thm:KS} in Section \ref{s:KS} below.

Consider the Szeg\H o recursion with the random Verblunski coefficients $\alpha_k$ from Theorem \ref{thm:KN}.  Proposition \ref{prop:unitary_car} yields the following result.
\begin{proposition}\label{prop:circ_car} Let $\beta>0$ and let $n$ be a positive integer.
Suppose that $\alpha_0, \dots, \alpha_{n-1}$ are distributed as in Theorem \ref{thm:KN},  and define $b_k, 0\le k\le n$, $b(t), t\in [0,1]$ according to  (\ref{car_center}).
The spectrum of the self-adjoint  Dirac operator (\ref{unop}) with  boundary conditions $u_0=(1,1)^t$ and $u_1=(b(1),1)^t$ is exactly the set $n \Lambda_n+2\pi n \Z$.
\end{proposition}

\bigskip

The following generalization of the circular ensemble appeared in \cite{ForBook} (see also \cite{FW2000}). Let $\beta>0$ and $\delta \in \CC$ with $\Re \delta >-1/2$. The finite circular Jacobi $\beta$-ensemble with parameter $\delta$ is a finite point process with joint probability density
\begin{equation}\label{huapickrell}
\frac{1}{Z_{n,\beta, \delta}} \prod_{1\le j<k\le n} |e^{i \lambda_j}-e^{i \lambda_k}|^\beta \prod_{j=1}^n \left|(1-e^{i\lambda_j})^{\overline \delta}\right|^2,\qquad \lambda_k\in [-\pi,\pi).
\end{equation}
For $\delta=k\beta/2$ with $k$ positive integer this can be viewed as the circular ensemble conditioned to have $k$ points at $1$. For $\beta=2$ these models were studied by \cite{Hua} and \cite{Pickrell}.

In \cite{BNR2009} the authors gave a construction for a random unitary matrix with joint eigenvalue density given by (\ref{huapickrell}) using the deformed Verblunsky coefficients.
\begin{theorem}[\cite{BNR2009}]\label{thm:BNR}
Let $\gamma_0, \dots, \gamma_{n-1}$ be independent random variables where the density of $\gamma_k$ is given by  $$c_{n,k} (1-|z|^2)^{\tfrac{\beta}{2}(n-k-1)-1}\left|(1-z)^{\bar \delta}\right|^2$$ on the unit disk for $k<n-1$, and $\gamma_{n-1}$ has density $$c_{n,n-1}\left|(1-z)^{\bar \delta}\right|^2$$ on the unit circle. Consider the $n\times n$ CMV matrix whose deformed Verblunsky coefficients are given by  $\gamma_0, \dots, \gamma_{n-1}$. Then the joint eigenvalue density of this matrix is given by (\ref{huapickrell}).
\end{theorem}
This is a generalization of the Killip-Nenciu construction. In the $\delta=0$ case the joint distribution of $\gamma_0, \dots, \gamma_{n-1}$ is the same as the one given in Theorem \ref{thm:KN} for the random Verblunsky coefficients. Moreover, because $|\alpha_k|=|\gamma_k|$ and the $\gamma_k$ random variables are rotationally invariant and independent, the `regular' Verblunsky coefficients $\alpha_k$ corresponding to the deformed coefficients $\gamma_0, \dots, \gamma_{n-1}$ will also have the same distribution as the one given in Theorem \ref{thm:KN}.

Proposition \ref{prop:unitary_car2b} yields a random Dirac operator where the spectrum is given by the periodic version of the circular Jacobi $\beta$-ensemble.

\begin{proposition}\label{prop:Hua}
Let $\beta>0$, $n$ a positive integer, and $\delta\in \CC$ with $\Re \delta>-1/2$.
Let $\gamma_0, \dots, \gamma_{k-1}$ be random variables with a joint distribution described as in Theorem \ref{thm:BNR}. Define $x_k, y_k$ via the recursion (\ref{xyrec}) with initial condition $x_0=0, y_0=1$ and set  $x(t)+i y(t)=x_{\lfloor nt \rfloor}+i y_{\lfloor nt \rfloor}$ for $t\in [0,1]$.
Then spectrum of the self-adjoint  Dirac operator  $\Dirop(x+iy, \infty, x(1))$  is exactly the set $n \Lambda_n+2\pi n \Z$.
\end{proposition}
 By the discussion around equation (\ref{xyrec}) we see that the matrices $X_k=\mat{1}{-x_k}{0}{y_k}$ satisfy the recursion $X_{k+1}=W_k X_k$ where $W_k$ are independent random matrices given by (\ref{AH}). Thus $X_k$ is a {\bf right random walk} on the affine group \eqref{affine}, and the matrix valued function $R_t$ of the corresponding Dirac operator is given by $R_t=\frac{1}{2\det X_k} X_k^t X_k$, with $k=\lfloor n t\rfloor$.

\subsection*{A special driving path}
The driving path for the circular $\beta$-ensemble has a simple intrinsic hyperbolic description. In this case, conditionally on the first $k$ steps, the point $x_{k+1}+iy_{k+1}$ is randomly chosen from uniform measure on the hyperbolic circle of random radius $r_k$ about the point $x_k+iy_k$.

This is also true for $k=n-1$, when $r_{n-1}=\infty$, $y_n=0$ and  $x_n$ is a point on the boundary chosen from harmonic measure centered at the point $x_{n-1}+iy_{n-1}$.

To see why, note that the increment of the walk comes from the rotationally invariant random variable $\alpha_k$ in the Poincar\'e disk representation, see \eqref{b_k}. There $|\alpha_k|^2$ has  Beta$(1,\tfrac{\beta}{2}(n-k-1))$ distribution, and $r_k=\log(1+|\alpha_k|)-\log(1-|\alpha_k|)$ is the the hyperbolic distance of $0$ and $\alpha_k$ in the Poincar\'e model.

\section{The Brownian carousel operator}  \label{s:bcar}

\cite{KS} showed that if $\Lambda_{n,\beta}$ is the finite point process with joint distribution (\ref{circdens})  then $n \Lambda_{n,\beta}$ has a point process scaling limit, the distribution of which can be described with a coupled system of SDEs.

Proposition \ref{prop:circ_car} shows that the set $n \Lambda_{n,\beta}+2\pi n\mathbb Z$ can be obtained  from a hyperbolic carousel driven by a random walk on the hyperbolic plane, or as the spectrum of the corresponding Dirac operator. The steps of this random walk are rotationally invariant, independent, and in the $n\to \infty$ limit the random walk path converges to a time-changed hyperbolic Brownian motion. This suggests that the scaling limit of $n \Lambda_{n,\beta}$ can be obtained from a hyperbolic carousel driven by a time changed hyperbolic Brownian motion and as the spectrum of the corresponding Dirac operator. In the upcoming work \cite{BVBV_future} we will  give a direct proof of this  using the convergence of the corresponding Dirac operators.  In this section we construct a random Dirac operator for the $\Sineb$ process, the scaling limit of the Gaussian $\beta$-ensemble. Because this process has already been described  via a hyperbolic carousel in \cite{BVBV}, the existence of the operator will follow from our results connecting Dirac operators and hyperbolic carousels.

\subsection{Gaussian $\beta$-ensemble and the $\Sineb$ process}
The Gaussian orthogonal and unitary ensembles are some of the most studied finite random matrix models. They are given as $\tfrac{1}{\sqrt{2}}(A+A^*)$ where $A$ is $n\times n$ matrix whose entries are i.i.d.~standard real (or complex) normal random variables. The resulting matrix has $n$ real eigenvalues with joint density  given by
\begin{align}\label{betadens}
\frac{1}{Z_{n,\beta}} \prod_{1\le i<j\le n} |\lambda_i-\lambda_j|^\beta e^{-\tfrac{\beta}{4}\sum_{i=1}^n \lambda_i^2},
\end{align}
where $\beta=1$ for the real and $\beta=2$ for the complex  case. The joint density (\ref{betadens}) makes sense for any $\beta>0$, the resulting distribution is called the Gaussian $\beta$-ensemble. Besides the classical $\beta=1, 2$ and $4$ cases (the last of which corresponds to real quaternion normals) there is no known invariant matrix model for these distributions. However one can construct a random tridiagonal model with the appropriate joint eigenvalue distribution, this is the result of \cite{DE}.

As $n\to \infty$ the support of (\ref{betadens}) is asymptotically $[-2\sqrt{n}, 2\sqrt{n}]$. In fact, if one rescales the empirical spectral measure by $n^{-1/2}$ then almost surely there exists a weak limit, the famous Wigner semicircle distribution with density $\rho(x)=\frac{1}{2\pi}\sqrt{(4-x^2)_+}$ (see e.g.~\cite{ForBook}). Let $|E|<2$ be a reference point inside the support of the limiting law. Then scaling the finite $\beta$-ensemble near $\sqrt{n} E$ by a factor of $\sqrt{4-E^2} \sqrt{n}$ one expects a limiting point process with asymptotic density $\tfrac{1}{2\pi}$. This has been proved to be true in the classical cases by Gaudin, Mehta and Dyson (see the monographs \cite{AGZ}, \cite{ForBook}, \cite{mehta}). For the general $\beta$ case we have the following.

\begin{theorem}[\cite{BVBV}]\label{thm:car}
Fix $\beta>0$ and $|E|<2$. Let $\Lambda_n$ be a finite point process with density (\ref{betadens}). Then $\sqrt{4-E^2} \sqrt{n} (\Lambda_n-\sqrt{n} E)$ converges in distribution to a point process $\Sineb$.
\end{theorem}

The  $\Sineb$ process has translation invariant distribution with intensity $\tfrac{1}{2\pi}$. As  the next theorem shows, the process $\Sineb$ can be obtained as the result of a hyperbolic carousel driven by a time changed hyperbolic Brownian motion. Recall the definition and basic properties of the hyperbolic Brownian motion from Subsection \ref{subs:hypB}.

For notational convenience we introduce the logarithmic time change function
\begin{align}\label{timech}
\tch(t)=-\log(1-t).
\end{align}

Let $\mathcal B_t$ be hyperbolic Brownian motion with variance $\frac{4}{\beta}$, and let $\eta_0$ be a fixed deterministic boundary point of $\HH$.   Consider the hyperbolic carousel with driving path is $\mathcal B_{\tch(t)}, t\in [0,1)$ and the starting point is $\eta_0$. For any fixed $\lambda\in \R$ let $\alpha_\lambda(t)$ denote the continuous lifting of the hyperbolic angle of $\eta_0, \mathcal B_{\tch(t)}$ and $r_\lambda(t)$ to $\mathbb R$ with $\alpha_\lambda(0)= 0$.

In \cite{BVBV} it was shown that for all $\lambda$ a.s.,
$$\alpha_\lambda(1)=\lim_{t\to 1} \alpha_\lambda(t)
$$
exists and is in $2\pi \mathbb Z$.

\begin{theorem}[The Brownian carousel, \cite{BVBV}]\label{thm:bcar}
The right continuous version of the total winding number  $\alpha_\lambda(1)/(2\pi)$ is the counting function of the $\Sineb$ process.
\end{theorem}

Note that one has to apply a simple time-change to get the actual form of the corresponding theorem in \cite{BVBV}.

\subsection{The $\Sineop$ operator}

In this subsection we first construct the Dirac operator $\Sineop$. Its oscillation theory is similar to the carousel representation of Theorem \ref{thm:bcar}.
 In the second part, we show that the spectum of the operator is indeed the $\Sineb$ process.

\begin{theorem}\label{thm:SSO}
Fix $\beta>0$.
Let $x+i y$ be the hyperbolic Brownian motion with variance $\frac{4}{\beta}$ in the half plane with initial condition $i$, as defined in (\ref{hypBM01}). Let $q=\lim_{t\to \infty} x(t)$, and use the notation $\tilde x(t)=x(\tch(t))$, $\tilde y(t)=y(\tch(t))$ for the time-changed process, see \eqref{timech}. Then the operator
\[
\Sineop=\Dirop(\tilde x+i \tilde y, \infty, q)
\]
on the interval $[0,1)$ satisfies conditions (\ref{cond1})-(\ref{cond3}), and hence it is self-adjoint on the appropriate domain.  The operator is limit circle for $\beta>2$, limit point for $\beta\le 2$. The inverse of the operator is a.s.~Hilbert-Schmidt.
\end{theorem}

\begin{proof}%
By (\ref{Rrep}) we have $$R=\frac{1}{2\tilde y}\mat{1}{-\tilde x}{-\tilde x}{\tilde x^2+\tilde y^2}.$$
The conditions (\ref{cond1}) and (\ref{cond2}) are satisfied since $\tilde x+i \tilde y$ is locally bounded in $\HH$.

Set $u_0=(-1,0)^t$, $u_1=(1,q^{-1})^t$.
Then we have
\begin{align}\label{uRu}
\int_0^1 u_1^t R(s) u_1 ds=\int_0^1 \frac{(q-\tilde x(s))^2+\tilde y(s)^2}{2 \tilde y(s)} ds=\frac{1}{2} \int_0^\infty  e^{-s} \left(\frac{(q-x(s))^2}{y(s)}+y(s)\right) ds.
\end{align}
We can write $y(t)=e^{\sigma W_1(t)-\sigma^2 t/2}$ and $x(t)=\sigma \int_0^t e^{\sigma W_1(s)-\sigma^2 s/2} dW_2(s)$ with independent standard Brownian motions $W_1, W_2$ and $\sigma^2=\frac{4}{\beta}$. The process $q-x(t)$ has the same distribution as $\sigma W(\int_t^\infty e^{2\sigma W_1(s)-\sigma^2 s} ds)$ where $W(\cdot)$ is a standard Brownian motion independent of $W_1$. For every $\delta>0$ standard Brownian motion satisfies
$$
|B(t)|\le C t^{1/2+\delta}
$$
for some random constant $C$ and all $t>1$.
It follows that for any $\eps>0$ there is a random constant $C$ so that the following inequalities hold with probability one:
\begin{align}\label{BM_asympt}
C^{-1} e^{-(\frac{2}{\beta}+\eps) t}\le y(t) \le C e^{-(\frac{2}{\beta}-\eps)t}, \qquad |q-x(t)|\le C e^{-(\frac{2}{\beta}-\eps) t}.
\end{align}
Choosing $\eps$ small enough we get that $\int_0^1 u_1^t R(s) u_1 ds<\infty$ a.s. This gives condition (\ref{cond3}) with $u_*=u_1$.

To see whether the Dirac operator in question is limit point or limit circle we need to check whether $\int_0^1 u^t R(s) u ds<\infty$ for a nonzero vector $u$ not parallel to $u_1$. Taking $u=(1,0)^t$ we get
\[
\int_0^1 u^t R(s) u ds=\frac{1}{2}\int_0^1 \tilde y(s)^{-1} ds=\frac12  \int_0^\infty e^{-s} y(s)^{-1} ds=\frac12  \int_0^\infty e^{-\frac{2}{\sqrt{\beta}} W_1(s)+(\frac{2}{\beta}-1)s} ds,
\]
which is a.s.~finite for $\beta>2$ (this is the limit circle case) and a.s.~infinite for $\beta\le 2$ (this is the limit point case).

To show that the inverse of the operator is Hilbert-Schmidt  we need to check that the integral \eqref{eq:HSnorm} is finite. Using the bounds \eqref{BM_asympt} we get
\begin{align*}
 \int_0^1 \int_0^t u_0^t R(s) u_0 \, u_1^t R(t) u_1 \,ds\, dt&=\frac14 \int_0^\infty \int_0^t e^{-(s+t)}\frac{1}{y(s)}  \left(\frac{(q-x(t))^2}{y(t)}+y(t)\right) ds\, dt\\
 &\le C' \int_0^\infty \int_0^t e^{-(s+t)+(\frac{2}{\beta}+\eps) s-(\frac{2}{\beta}-3\eps)t} ds\, dt,
\end{align*}
which is finite if $\eps>0$ is small enough. This completes the proof of the theorem.
\end{proof}

Theorem \ref{thm:SSO} describes the $\Sineop$ operator via the hyperbolic Brownian motion in the upper half plane. One can also use the disk representation of the Brownian motion.  For this we consider the similar operator $\tilde U \circ \Sineop\circ \tilde U^{-1}$ where $\tilde U=\tfrac{1}{\sqrt{2}} \mat{1}{-i}{1}{i}$ is the  linear transformation corresponding to the Cayley map between  the unit disk and the half plane representation. The resulting operator is
\[
\widetilde{\Sineop}  u=\frac{2}{(1-|b(t)|^2)} \mat{1+|b(t)|^2}{2 b(t)}{2 \bar b(t)}{1+|b(t)|^2} \mat{-i}{0}{0}{i} u',
\]
where $b$ is the hyperbolic Brownian motion with variance $\frac{4}{\beta}$ in the Poincar\'e disk after the time change $\tch(t)$. Note that $u$ is a $\CC^2$ valued function on $[0,1)$.

\begin{figure}[ht]\hskip20mm
\includegraphics*[width=330pt]{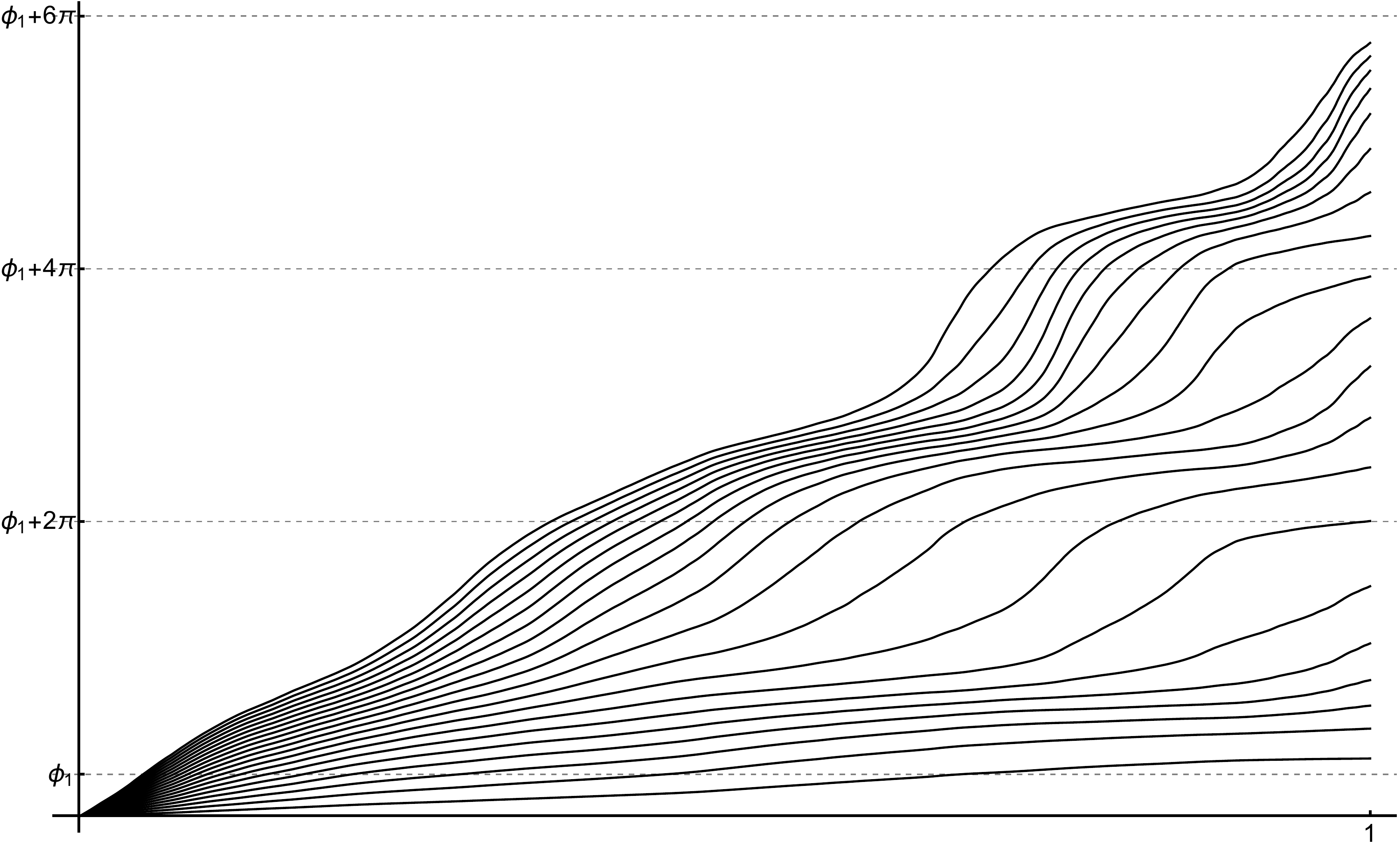}\\
\caption{Simulation of the phase function of $\Sineop$ with $\beta>2$ with various values of $\lambda$}
\end{figure}

\begin{figure}[ht]\hskip20mm
\includegraphics*[width=330pt]{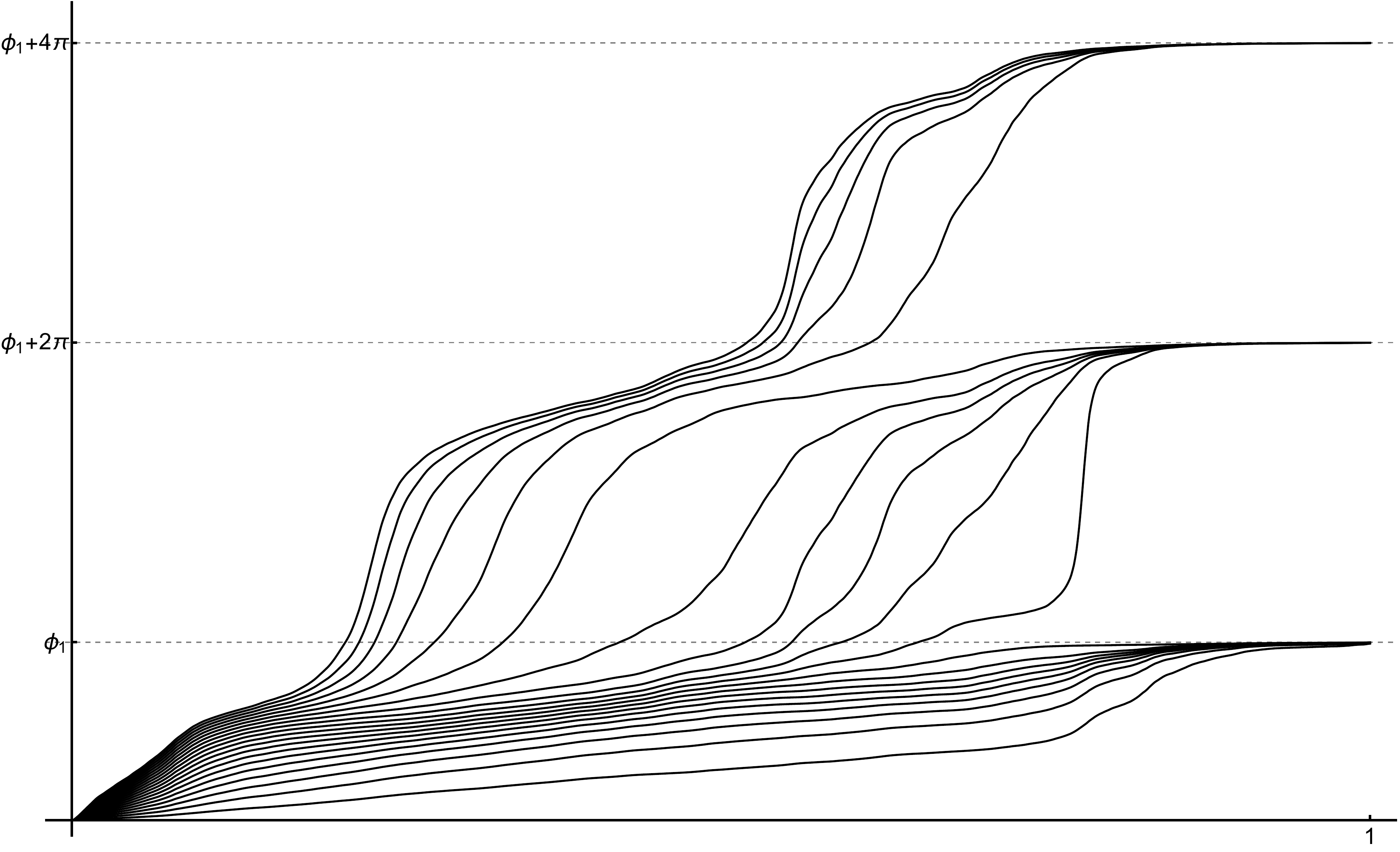}\\
\caption{Simulation of the phase function of $\Sineop$ with $\beta\le2$ with various values of $\lambda$}
\end{figure}

\begin{theorem}\label{thm:SSO1}
Fix $\beta>0$ and consider the $\Sineop$ operator defined in Theorem \ref{thm:SSO}. The spectrum of $\Sineop$ is distributed as the $\Sineb$ process
\end{theorem}

\begin{proof}
By Theorem \ref{thm:SSO} the spectrum of the $\Sineop$ operator is a.s.~a simple point process. By the oscillation theory in Section \ref{subs:phase}, the counting function is the right-continuous version of the function $\lambda \to \textup{sign}(\lambda)\cdot  \left|\{t\in (0,1): r_\lambda(t)=q\}\right|$ where $r_\lambda(t)$ is the moving boundary point in the carousel given by the ODE system
\[
r_\lambda'(t)=\lambda   \frac{\tilde y(t)^2+(r_\lambda(t)-\tilde x(t))^2}{2 \tilde y(t)}, \qquad r_\lambda(0)=-\infty.
\]
Recall Theorem \ref{thm:bcar} and the notation there. Let $\tilde x(t)+i\tilde y(t)$ be $B_{\tch(t)}$ in half-plane coordinates. Then the hyperbolic angle $\alpha_\lambda$ between the points $\infty, \cB_{\tch(t)}, r_\lambda(t)$ satisfies
\begin{align}\label{r_alpha}
\cot \left(\frac{\alpha_\lambda(t)}{2}\right)=\frac{\tilde x(t)-r_\lambda(t)}{\tilde y(t)}.
\end{align}
It suffices to show that for all $\lambda\in \R$ a.s.
\begin{align}\label{counting}
\textup{sign}(\lambda) \cdot \left|\{0<t<1: r_\lambda(t)=q\}\right|=\frac{\alpha_\lambda(1)}{2\pi},
\end{align}
this implies that the counting function of the spectrum of the $\Sineop$ operator is given by the $\Sineb$ process.

We show (\ref{counting}) for $\lambda>0$, the other case follows similarly. Moreover, we will consider the $\beta\le 2$ and $\beta>2$ cases separately.

If $\beta\le 2$ then the  $\Sineop$  operator is limit point at $t=1$, so by Theorem \ref{thm:osc} we have  $\lim_{t\to 1} r_\lambda(t)=q$ for all $t\neq 0$. Assume $\lambda>0$, the other case will follow similarly. Then $r_\lambda(t)$ is strictly increasing in $t$  (between blowups), and is continuous as a function to $\partial \mathbb H=\mathbb R\cup \{\infty\}$. This implies that $r_\lambda(t)$ converges to $q$ from below. As a consequence, the following cardinalities are equal:
$$\left|\{0<t<1: r_\lambda(t)=q\}\right|=\left|\{0<t<1: r_\lambda(t)=\infty\}\right|.$$ By (\ref{r_alpha}) the blowup times of $r_\lambda(t)$ are exactly the times when $\alpha_\lambda(t)$ hits an integer multiple of $2\pi$. The function $\alpha_\lambda(t)$ converges to an integer multiple of $2\pi$, and by Theorem 7 of \cite{BVBV} for $\beta\le 2$ the function $t\to \alpha_\lambda(t)$ will converge to its limit from above. By Proposition 9 of \cite{BVBV} the process $\alpha_\lambda(t)$ cannot go below $2\pi n$  ($n\in \Z$) once it hits this value, which means that  $\lim_{t\to 1}\frac{1}{2\pi}\alpha_\lambda(t)$ is exactly the same as the the number of blowups of $r_\lambda(t)$. This completes the proof of the proposition for the $\beta\le 2$ case.

For $\beta>2$ the operator $\Sineop$ is limit circle at $t=1$. Theorem \ref{thm:osc} implies that $r_\lambda(1)=\lim_{t\to 1} r_\lambda(t)$ exists and the limiting function $\lambda\to r_\lambda(1)$ is continuous and
 strictly increasing   (apart from blowups at $\infty$). If $\lambda>0$ is not an eigenvalue then $r_\lambda(1)=\lim_{t\to 1} r_\lambda(t)$ is not equal to $q$. If $r_\lambda(1)<q$ then $t\to r_\lambda(t)$ hits $q$ the same number of times as it hits $\infty$.  Equation (\ref{r_alpha}) shows that $\alpha_\lambda(t)$ converges to an integer multiple of $2\pi$ from above which means that the number of blowups of $r_\lambda$ is the same as the limit of $\frac{1}{2\pi}\alpha_\lambda(t)$ as $t\to1$.
If  $r_\lambda(1)>q$ then the size of the set $\{0<t<1: r_\lambda(t)=q\}$ is equal to the number of blowups of $r_\lambda(t)$ plus one. But in that case  $\alpha_\lambda(t)$ converges to an integer multiple of $2\pi$ from below (by equation (\ref{r_alpha})), which means that the number of blowups of $r_\lambda(t)$ is equal to $\lim_{t\to 1} \frac{1}{2\pi}\alpha_\lambda(t)-1$. This shows that if $\lambda$ is not an eigenvalue then the size of the set $\{0<t<1: r_\lambda(t)=q\}$ is  equal to $\lim_{t\to 1} \frac{1}{2\pi}\alpha_\lambda(t)$, which finishes the proof.
\end{proof}

\section{The Killip-Stoiciu limit of the circular $\beta$-ensemble is the $\Sineb$ process} \label{s:KS}

\cite{KS} show that the circular $\beta$-ensemble has a point process scaling limit which can be characterized by a coupled system of SDEs.
\begin{theorem}[\cite{KS}]\label{thm:KS}
Fix $\beta>0$ and let $\Lambda_n$ be the finite point process with density given by (\ref{circdens}). Then $n \Lambda_n$ converges in distribution to the point process
 $$\Xi = \{\lambda: \psi_\lambda(1)\in  \theta+ 2\pi\Z\}$$
where $\psi_\lambda(t)$ is the strong solution of the one-parameter family of SDEs
\begin{align}\label{KillipSDE}
d\psi_\lambda=\lambda dt+\frac{2}{\sqrt{\beta t}} \Re\left[ (e^{-i \psi_\lambda}-1)(dB_1+i dB_2)\right], \qquad \psi_\lambda(0)=0, \quad t\in [0,1],
\end{align}
with $\lambda \psi_\lambda(t)\ge 0$ for all $\lambda, t$.
Here $B_1, B_2$ are standard Brownian motions, $\theta$ is uniform on $[0,2\pi]$ and the three are independent.
\end{theorem}
The SDE system (\ref{KillipSDE}) looks  similar to the system used to describe the winding angle $\alpha$ of Theorem \ref{thm:bcar} about the $\Sineb$ process. In  \cite{BVBV} it was shown that $\alpha$ satisfies
\begin{align}\label{sse}
d\alpha_\lambda=\lambda dt+\frac{1}{\sqrt{1-t}} \frac{2}{\sqrt{\beta}} \Re \left[(e^{-i \alpha_\lambda}-1)(dB_1+i dB_2) \right], \qquad \alpha_\lambda(0)=0,
\end{align}
The SDE \eqref{sse} resembles a time-reversed version of \eqref{KillipSDE}.

For the classical values $\beta=1, 2$ and 4 the point process limit of the circular beta ensemble is the same as the $\Sineb$ process, see \cite{ForBook}. In \cite{Nakano} this was shown for general $\beta>0$ by deriving both processes as scaling limits of the spectra of certain Schr\"odinger operators.

The first item on the list of open problems from the 
\cite{AIM} is to describe the direct connection between the two characterizations.
We obtain both SDE systems from the same $\Sineop$ operator by considering the ordinary and reverse oscillation theories. In this coupling the associated point processes satisfy $\Sineb=-\Xi$. Since $\Sineb\ed-\Sineb$, this coupling shows that $\Xi\ed\Sineb$.

Consider the $\Sineop$ operator $\Dirop(\mathcal B, \infty, \cB(1))$, associated to the hyperbolic Brownian motion run in logarithmic time, as constructed in Theorem \ref{thm:SSO}.

Let $\phi_t$ be the phase function in the oscillation theory, and let $\rho_t$ be the reverse phase function, see Section \ref{s:carousel}.
Recall that
$$
\alpha_{\lambda}(t)=\rm{angle}(\infty, \cB(t),-\cot(\phi_\lambda(t)/2)),\qquad\alpha_\lambda(0)=0,
$$
and that $\alpha_\lambda(1)/(2\pi)$ is the counting function of
$\spec \Sineop$. Let
$$
\psi_{\lambda}(t)=\rm{angle}(\cB(1),\cB(1-t),-\cot(\rho_{-\lambda}(1-t)/2))\qquad\psi_\lambda(0)=0,
$$
be the continuous lifting of this hyperbolic angle.
\renewcommand{\theenumi}{(\roman{enumi})}
\renewcommand{\labelenumi}{\theenumi}
\begin{theorem}[The Killip-Stoiciu SDE and the Sine$_\beta$ process] Let $\theta=\rm{angle} (\cB(1),i,\infty)$.
\begin{enumerate}
\item A.s.~we have
$\spec \Sineop = \{-\lambda: \psi_{\lambda}(1)= \theta \,\,\, \textup{mod } 2\pi\}.$
\item The process $\psi_\lambda(t)$ and the angle $\theta$ have joint distribution described in Theorem \ref{thm:KS}.
    \end{enumerate}
\end{theorem}

\begin{proof}
%
For (i) note that by reverse oscillation theory and by the definition of $\psi$ we have
$$
\spec \Sineb= \{\lambda: \rho_{\lambda}(0)\,\in \,2\pi \mathbb Z\}=\{-\lambda: \psi_{\lambda}(1)\,\in\, \theta +2\pi \mathbb Z\}.
$$
The existence and strong uniqueness for the SDE system (\ref{KillipSDE}) is
proved in \cite{KS}. The standard theory does not apply because of the blowup of the  diffusion coefficient at 0. However it is possible to approximate this system by `nicer'  systems. We outline  a version of the argument in \cite{KS}. Consider the SDE system
\begin{align}\label{KillipSDE2}
d\psi_{\lambda,\eps}=\lambda dt+\frac{2}{\sqrt{\beta t}} \Re\left[ (e^{-i \psi_{\lambda,\eps}}-1)(dB_1+i dB_2)\right], \qquad \psi_{\lambda,\eps}(\eps)=0,\quad t\in [\eps,1].
\end{align}
This system has a unique strong solution $\psi_{\lambda,\eps}(t)$ which is a continuous function of both $\lambda$ and $t$
as the coefficients are globally Lipschitz continuous.

For any $\lambda\neq 0$ we have $\lambda \psi_{\lambda,\eps}(t)>0$ for $t>\eps$ and if $\eps_1>\eps_2$ then $\lambda(\psi_{\lambda,\eps_2}(t)-\psi_{\lambda,\eps_1}(t))>0$ for $t>\eps_1$. These statements follow from the fact that two solutions of the SDE (\ref{KillipSDE}) (ignoring the initial condition) that are ordered at a certain time $t_0$, are also ordered the same way for $t\ge t_0$. Extend the definition of $\psi_{\lambda,\eps}(t)$ for the full $t\in[0,1]$ interval by defining the process to be 0 on $[0,\eps]$.  Then the monotone limit of the solutions $\psi_{\lambda,\eps}$ as $\eps\to 0$ gives the strong solution of (\ref{KillipSDE}) on $[0,1]$.

Recall that $\rho_{-\lambda}$ is the monotone limit of $\rho_{-\lambda,1-\eps}$ as $\eps \to 0$. Let
$$
\psi_{\lambda,\eps}(t)=\rm{angle}(-\cot(\rho_{-\lambda,1-\eps}(1-t)/2),\cB(1-t),\cB(1))\qquad\psi_{\lambda,\eps}(\eps)=0,
$$
then as $\eps\to 0$ the function $\psi_{\cdot, \eps}$ converges to $\psi_\cdot$ given in \eqref{KillipSDE}. We will show that $\theta$ is independent of $\psi_{\cdot,\eps}$ and that the latter satisfies the SDE \eqref{KillipSDE2} for some Brownian motions $B_1,B_2$. By the argument above this shows that $\psi_{\cdot, \eps}$ converges in distribution to the solution of \eqref{KillipSDE}. Together with the independence of $\theta$ this implies (ii).

Note that $\rho_{\cdot , 1-\eps}$ is the reverse phase function for $\Dirop(\mathcal B(t), t\in [0,1-\eps], \infty, \cB(1))$.
Let $\MT_\eps$ be the M\"obius transformation taking
$\cB(1-\eps)\mapsto i$ and $\cB(1)\mapsto \infty$. We reverse time in the above operator and apply $\MT_\eps$ to its components. Let $\xi$ be the phase function for the resulting operator $\Dirop(\MT_\eps\mathcal B(1-t), t\in [\eps,1], \MT_\eps\cB(1),\MT_\eps \infty)$. The time interval starts from $\eps$ and not the usual $0$. Recall that $\psi_{\lambda,\eps}$ is the hyperbolic angle of the moving boundary point of the carousel, the center of rotation, and reference point. Its evolution is invariant under hyperbolic isometries. In the time-reversed operator it becomes the ordinary phase function, and $-\lambda$ gets replaced by $\lambda$.  This gives
$$
\psi_{\lambda,\eps}(t)=
\rm{angle}(-\cot(\xi_{\lambda,\eps}(t)/2),\MT_\eps\cB(1-t),\MT_\eps\cB(1)).
$$
By Proposition \ref{prop:reverse} the law of $(\MT_\eps\cB(1-t), t\in (\eps,1))$ is hyperbolic Brownian motion with the appropriate local variance $\tfrac{4}{\beta}\,\tch(1-t)$. Moreover, by the Propostion it is independent of the uniform angle $(\infty, \MT_\eps \cB(0), \MT_\eps \infty)$. This a agrees with the same angle before the transformation: $(\cB(1),i,\infty)=\theta$. This implies that $\psi_{\lambda,\eps}$ satisfies the SDE analogous to the one satisfied by $\alpha$, namely the SDE \eqref{KillipSDE2}.
\end{proof}

\section{Classification of operator limits of random matrices  }
\label{s:zoo}

In this section we consider a family of stochastic Dirac operators generalizing the Brownian Carousel operator $\Sineop$. We start by considering specific examples connected to various random matrix models and then we discuss the general family.

\subsection{The hard edge operator}

The Laguerre $\beta$-ensemble is a generalization of the gaussian Wishart matrices. It is a two-parameter family of finite ensembles given by the following density function:
\begin{align}\label{lagdens}
\frac{1}{Z_{n,\beta,a}} \prod_{1\le i<j\le n} |\lambda_i-\lambda_j|^\beta \prod_{j=1}^n \lambda_j^{\tfrac{\beta}{2}(a+1)-1}e^{-\tfrac{\beta}{2} \lambda_j}, \qquad \lambda_j\ge 0.
\end{align}
Here $n$ is a positive integer and $a>-1$. For $\beta=1,2,4$ and integer $a$ one can realize this ensemble as the eigenvalues of a matrix $MM^*$ where $M$ is $n\times (n+a)$ with i.i.d.~standard gaussian entries (with real, complex or real quaternion random variables).

If $a>-1$ is kept fixed and $n\to \infty$, the support of the finite point process will be asymptotically $[0,4n]$ and the limiting empirical spectral density (after rescaling by $\frac{1}{n}$) will be given by  $\frac{1}{2\pi}\left(\frac{4-x}{x}\right)^{1/2}\cdot 1_{\{x>0\}}$. The limit of the process scaled by $n$ is expected to be different from the already discussed bulk  limit. In the classical $\beta=1,2$ and 4 cases the limits where derived and characterized by Tracy and Widom \cite{TW_hard}. The general $\beta$ case is the following.
\begin{theorem}[Hard edge limit, \cite{RR}, \cite{RR_err}]\label{thm:RR}
Fix $a>-1$ and  $\beta>0$. Let $\Lambda_{n}$ be the finite non-negative point process with joint density (\ref{lagdens}). Then $n \Lambda_n$ converges to a simple point process, namely the discrete spectrum of the following Sturm-Liouville differential operator:
\begin{align}\label{hardop}
\mathfrak{G}_{\beta,a} f(x)= -\frac{1}{m(x)} \partial_x \left(\frac{1}{s(x)} \partial_x \, f(x) \right), 
\end{align}
acting on functions $[0,\infty)\to \R$ with Dirichlet boundary condition at 0  and Neumann boundary condition at $\infty$. Here $B(x)$ is a standard Brownian motion and
\begin{align}
m(x)=e^{-(a+1)x-\frac{2}{\sqrt{\beta}}B(x)}, \quad s(x)=e^{ax+\frac{2}{\sqrt{\beta}}B(x)}.
\end{align}
\end{theorem}
The theory of Sturm-Liouville operators is closely connected to that of the Dirac operators. (See \cite{Weidmann} or Chapter 9 of \cite{Teschl}.) The operator $\mathfrak{G}_{\beta,a}$ is self-adjoint with domain given by the following subset of $L^2_m=L^2(\R_+, m\, dx)$:
\begin{align}\label{domG}
\dom_{\mathfrak{G}_{\beta,a}}=\{f:\R^+\to \R: f, s^{-1} f'\in \ac(\R^+), f, \mathfrak{G}_{\beta,a}f\in L^2_m, f(0)=0, \lim_{x\to \infty} s^{-1}(x)f'(x)=0   \}.
\end{align}
The inverse of the operator is a Hilbert-Schmidt integral operator defined as
\begin{align}\label{hardinverse}
\mathfrak{G}_{\beta,a}^{-1} f(x)=\int_0^\infty K(x,y) f(y) m(y) dy, \quad K(x,y)=\int_0^xs(z) dz \ind(x<y)+\int_0^ys(z) dz \ind(y\le x).
\end{align}
As we will show in the next theorem, the Sturm-Liouville operator $\mathfrak{G}_{\beta,a}$ can be transformed into a  Dirac operator that fits into our framework. The corresponding  hyperbolic carousel is driven by a real Brownian motion with a drift on a line in the hyperbolic plane.

The Euclidean real line $\R$ is embedded into the imaginary axes $\{iy: y>0\}\subset\HH$ by the transformation $x\to i e^x$. A real Brownian motion with drift in this embedded line is just geometric Brownian motion moving on the set $i \R_+$.

\begin{theorem}\label{thm:hardOP}
Fix $\beta>0$, $a>-1$ and let $B$ be standard Brownian motion. Let $y(t)=e^{\frac{2}{\sqrt \beta}B(2t)+\left(2a+1\right)t}$ and define $\tilde y(t)=y(\tch(t))$ with $\tch$ from (\ref{timech}). Then the operator
\[
\Besselop=\Dirop(i \tilde y, 0, \infty)
\]
on the interval $[0,1)$ satisfies the conditions (\ref{cond1})-(\ref{cond3}) and hence it is self-adjoint on the appropriate domain. The operator is  limit circle near 1  for $-1<a<0$ and limit point for $a\ge 0$. It is a.s.~invertible and the inverse is a.s.~Hilbert-Schmidt. Moreover, we have the equality of spectra $$\spec \Besselop=4\sqrt{\spec \mathfrak{G}_{\beta,a}}\cup (-4\sqrt{\spec \mathfrak{G}_{\beta,a}})$$ where $\mathfrak{G}_{\beta,a}$ is the operator \eqref{hardop} built from $B$.
\end{theorem}
\begin{proof}
Set
\begin{equation}\label{Qdef}Q(t)=e^{-\frac{1}{2} t}\mat{\frac{1}{y(t/2)}}{0}{0}{y(t/2)}=\mat{m(t)}{0}{0}{s(t)},
\end{equation}
and consider the Dirac operator $\kappa_{\beta,a}=Q^{-1} J \frac{d}{dx}$ on $[0,\infty)$. Note that since $a>-1$ we have
\begin{align}
\int_0^\infty (1,0)Q(x) (1,0)^t dx=\int_0^\infty m(x) dx=\int_0^\infty e^{-(a+1)x-\frac{2}{\sqrt{\beta}}B(x)}dx<\infty\quad \textup{a.s},
\end{align}
and $\int_0^\infty (0,1)Q(x) (0,1)^t dx=\int_0^\infty s(x) dx=\int_0^\infty e^{ax+\frac{2}{\sqrt{\beta}}B(x)}dx<\infty$ if and only if $-1<a<0$. Thus $\kappa_{\beta,a}$ satisfies conditions (\ref{cond1})-(\ref{cond3}), it is limit point at $\infty$ for $a\ge 0$ and limit circle there  for $-1<a<0$. The operator is self-adjoint on the appropriate domain with initial condition $u_0=(0,1)^t$ and end condition $u_1=(1,0)^t$. We also have
\begin{align*}
 \int_0^\infty \int_0^t u_0^t Q(s) u_0 \, u_1^t Q(t) u_1 ds dt=\int_0^\infty \int_0^t e^{-\frac{1}{2}(t+s)} e^{\frac{2}{\sqrt{\beta}}(B(s)-B(t))+(a+\frac12)(s-t)}ds dt,
 \end{align*}
which is a.s.~finite since $\lim_{t\to \infty}\frac{B(t)}{t}=0$ and $a>-1$. This shows that $\kappa_{\beta,a}$ is a.s.~invertible with a Hilbert-Schmidt inverse.

Using (\ref{hardinverse}) and the Cauchy-Shwarz inequality one can check that $f$ is in the domain of the operator $\mathfrak{G}_{\beta,a}$ if and only if $(f,0)^t$ and $(0,s^{-1} f')^t$ are in the domain of the operator ${\kappa_{\beta,a}}$. From the definition we get that if $\mathfrak{G}_{\beta,a} f=\nu^2 f$ then
\begin{align}\label{xxx123}
\kappa_{\beta,a} \bin{f_1}{f_{2}}=\pm \nu \bin{f_1}{f_2}, \qquad \textup{with} \qquad f_1=f, \quad f_{2} =\pm \nu^{-1} s^{-1} f'.
\end{align}
Moreover, if $(f_1,f_2)^t$ satisfies $\kappa_{\beta,a} (f_1,f_2)^t=\nu (f_1,f_2)^t$ then $\mathfrak{G}_{\beta,a} f_1=\nu^2 f_1$.

This implies that the spectrum of $\kappa_{\beta,a}$ is $(-\sqrt{\spec \mathfrak{G}_{\beta,a}})\cup \sqrt{\spec \mathfrak{G}_{\beta,a}}$.
To make the connection more precise consider the following subspaces in the domain of $\kappa_{\beta,a} $:
\[
\mathcal H_{\pm}=\left\{ \left(\pm f, s^{-1} \left(\mathfrak{G}_{\beta,a} ^{-1/2} f\right)'\right)^t : f\in \dom_{\mathfrak{G}_{\beta,a}} \right\}.
\]
The computations around (\ref{xxx123}) show that $\kappa_{\beta,a}$ is isometric to $\mathfrak{G}_{\beta,a} ^{1/2}$ on $\mathcal H_{+}$  and isometric to  $-\mathfrak{G}_{\beta,a} ^{1/2}$ on $\mathcal H_{-}$. Moreover, $\mathcal H_{+}$ and $\mathcal H_{-}$ span the domain of $\kappa_{\beta,a}$.

To finish the proof we observe that the time-change $t\mapsto 2 \tch(t)=-2 \log(1-t)$ maps $\kappa_{\beta,a}$ to the Dirac operator $\tfrac{1}{4}\Besselop$ from which the theorem follows.
\end{proof}

\subsection{The Hua-Pickrell operator}

In Section \ref{s:circ} we considered the  finite circular Jacobi $\beta$-ensemble (\ref{huapickrell}), a  generalization of the circular $\beta$-ensemble. In Theorem \ref{prop:Hua} we showed that the point process can be obtained from a hyperbolic carousel driven by an affine random walk. By studying the asymptotic step distribution of the random walk one can show that the path converges in distribution to $x(\tch(t))+i y(\tch(t))$,  where $\tch(t)$ is given in (\ref{timech}) and  the process $x+i y$ solves the SDE
\begin{align}\label{Huasde}
dy=\left(-\frac{4}{\beta}\Re \delta dt+\frac{2}{\sqrt{\beta}}dB_1\right)y, \quad
dx=\left( \frac{4}{\beta}\Im \delta dt+\frac{2}{\sqrt{\beta}}dB_2\right)y, \qquad y(0)=1, x(0)=0.
\end{align}
This suggests that the appropriately scaled circular Jacobi $\beta$-ensemble converges in distribution to a point process which can be obtained by a random Dirac operator, or equivalently, it can be obtained from a hyperbolic carousel driven by  $x(\tch(t))+i y(\tch(t))$.  The rigorous proof of this statement will be given in the forthcoming paper \cite{BVBV_future}. Here we only give the description of the limiting Dirac operator.

\begin{proposition}\label{prop:Huaop}
Let $\beta>0$ and $\delta\in \CC$ with $\Re \delta>-1/2$.  Consider the solution $x, y$ of the SDE system (\ref{Huasde}) and let $\tilde x(t)=x(\tch(t))$, $\tilde y(t)=y(\tch(t))$ where $\tch$ is defined in (\ref{timech}).
Then the following statements hold:
\begin{enumerate}[(a)]
\item The limit $q=\lim_{t\to \infty} x(t)\in \R$ exits a.s.

\item The Dirac operator $\Huaop(\tilde x+i \tilde y, \infty, q)$ on $[0,1)$ satisfies conditions (\ref{cond1})-(\ref{cond3}), and hence it is self-adjoint on the appropriate domain. The operator is limit circle near 1 for $\Re \delta+\tfrac12<\tfrac{\beta}{4}$  and limit point for $\Re \delta+\tfrac12\ge \tfrac{\beta}{4}$.

\item The operator $\Huaop$ is a.s.~invertible and the inverse is a.s.~Hilbert-Schmidt.

\end{enumerate}
\end{proposition}
\begin{proof}
The SDE system  (\ref{Huasde}) can be solved explicitly to give
\[
y=e^{\frac{2}{\sqrt{\beta}}B_1(t)-\frac{4}{\beta}(\Re \delta+\tfrac12) t }, \qquad x=\frac{2}{\sqrt{\beta}}\int_0^t e^{\frac{2}{\sqrt{\beta}}B_1(s)-\frac{4}{\beta}(\Re \delta+\frac12) s } dB_2+\frac{4}{\beta}\Im \delta \int_0^t e^{\frac{2}{\sqrt{\beta}}B_1(s)-\frac{4}{\beta}(\Re \delta+\frac12) s } ds.
\]
Since $\Re \delta+\tfrac12>0$ we immediately get that $y(t)\to 0$ a.s.~as $t\to \infty$. We also get that the limit of $x(t)$ exists a.s.~and it is equal to
\[
q=\frac{2}{\sqrt{\beta}}\int_0^\infty e^{\frac{2}{\sqrt{\beta}}B_1(s)-\frac{4}{\beta}(\Re \delta+\frac12) s } dB_2+\frac{4}{\beta}\Im \delta \int_0^\infty e^{\frac{2}{\sqrt{\beta}}B_1(s)-\frac{4}{\beta}(\Re \delta+\frac12) s } ds.
\]
The rest of the proof follows the strategy of the proof of Theorem \ref{thm:SSO}. We can prove that for any small enough $\eps>0$ there is a random positive constant $C$ so that
\begin{align}
C^{-1} e^{-\frac{4}{\beta}(\Re \delta+\frac12+\eps) t}\le y(t) \le C e^{-\frac{4}{\beta}(\Re \delta+\frac12-\eps)t}, \qquad |q-x(t)|\le C e^{-\frac{4}{\beta}(\Re \delta+\frac12-\eps) t}.\label{xxxhua}
\end{align}
From this we get that
\[
\int_0^1 u_1^t R(s) u_1 ds=\frac{1}{2} \int_0^\infty e^{- s} \left(\frac{(q-x(s)^2}{y(s)}+y(s)\right) ds<\infty\qquad \textup{a.s.}
\]
We also have
\[
\int_0^1 u_0^t R(s) u_0 ds=\frac{1}{2} \int_0^\infty  e^{-s}e^{-\frac{2}{\sqrt{\beta}}B_1(s)+\frac{4}{\beta}(\Re \delta+\frac12) s }ds,
\]
which is finite a.s.~exactly if $\Re \delta+\tfrac12<\tfrac{\beta}{4}$. This completes the proof of (b). To prove (c) we need to show
$\int_0^1 \int_0^t u_0^t R(s) u_0 \, u_1^t R(t) u_1 ds dt<\infty$ a.s., which follows from the bounds (\ref{xxxhua}).
\end{proof}

\subsection{Discrete random Schr\"odinger operators}

Fix $\sigma>0$, and let $\omega_k$ be i.i.d.~standard normal random variables.
 Consider the $n\times n$ tridiagonal matrix $H_{n,\sigma}$, where the off-diagonal terms are constant 1 and the $k^{th}$ diagonal element is $\frac{\sigma}{\sqrt{n}}\omega_k$. When $\sigma=0$ the empirical spectral distribution of these matrices would converge to an arcsine law on $[-2,2]$.
The limiting bulk eigenvalue distribution was given terms of the  hyperbolic carousel.

\begin{theorem}[Theorem 6 of \cite{KVV}]
Let  $0<\alpha<\pi/2$, let $U$ be a uniform random variable on $[0,2\pi]$, independent of the $\omega_k$, and set $\nu=\frac{\sigma^2}{\sin^2\alpha}$.
Then we have
$$n \sin \alpha (\spec H_{n,\sigma}-2\cos \alpha)-U\Rightarrow \Sch^*_\nu$$ in distribution.
The point process $\nu^{-1} \Sch^*_\nu$ is given by the hyperbolic carousel $\mathcal{HC}(b(t),\eta_0, b(\infty))$ on the time interval $[0,\nu)$, where $b(t)$ is standard hyperbolic Brownian motion and $\eta_0\in \partial \HH$ is fixed. \end{theorem}
%

%

Let $b=x+i y$ be in the half-plane representation, and let  $q=b(\infty)=x(\infty)$.
Then the Dirac operator $\Schop =\Dirop(x+i y, \infty, q)$ on the interval $[0,\nu)$ satisfies the conditions (\ref{cond1})-(\ref{cond3}),  self-adjoint on the appropriately defined domain,  and  $\nu\cdot\spec \Schop=\Sch^*_\nu$.

\subsection{Brownian motion on the affine group and random operators}

We have seen that the $\Sineop$, $\Besselop$ and $\Huaop$ operators are all of the form (\ref{Diracop}).
 In all three cases the matrix valued stochastic process $R(t)$ is given as $\frac{1}{2 y}\mat{1}{- x}{- x}{ x^2+ y^2}$ under the logarithmic time change (\ref{timech}), where $x+i y$ is a certain diffusion on the upper half plane.  Moreover, the diffusions in question are all of the following form:
\begin{align}\label{gendiff}
dy=\left( \gamma_1 dt+ \alpha_1 dB_1\right)y, \quad
dx=\left( \gamma_2 dt+\alpha_2 dB_2\right)y, \qquad y(0)=1,\, x(0)=0.
\end{align}
Here $B_1$ and $B_2$ are independent standard Brownian motions.
The various values of these parameters are summarized in the table below.
\begin{align*}
\begin{array}{|c|c|c|c|c|}
\hline
& \alpha_1^2 &\alpha_2^2 & \gamma_1& \gamma_2\\
\hline
\Sineop & \frac{4}{\beta} \phantom{\Bigg|}& \frac{4}{{\beta}} &  0& 0\\
\hline
\Besselop & {\frac{8}{\beta}} \phantom{\Bigg|}& 0 & (2a+1)-\frac{4}{\beta}& 0\\
\hline
\Huaop & \frac{4}{{\beta}}& \frac{4}{{\beta}}  \phantom{\Bigg|}& -\frac{4}{\beta} \Re \delta &\frac{4}{\beta}\Im \delta\\
\hline
\end{array}
\end{align*}
The path $ x+i  y$ can be identified with a path $ X=\mat{1}{- x}{0}{ y}$ on the group of affine matrices.  $X$ is right Brownian motion on the group of affine matrices satisfying the SDE
\begin{align}\label{XXXSDE}
dX= d\mathcal B X, \qquad X_0=I, \qquad d\mathcal B=\mat{0}{ -\gamma_2 dt-\alpha_2 dB_2}{0}{\gamma_1 dt+ \alpha_1 dB_1}
\end{align}
The  matrix valued function $R(t)$ appearing in the operators is given by $\frac{ X^t X}{2\det  X}$ (under the logarithmic time change) which is just half the positive definite representation of the diffusion $x+iy$.

By It\^ o's classification of Brownian motion on a Lie group \cite{Ito50} any right Brownian motion on the group of affine matrices of the form $\mat{1}{-x}{0}{y}$ will satisfy the SDE  (\ref{XXXSDE}) with some choice of parameters, and possibly correlated standard Brownian motions $B_1$, $B_2$.  The Brownian motion corresponding to the $\Sineop$, $\Besselop$ and $\Huaop$ operators does not cover all possible parameter values. It would be interesting to see whether there exists random matrix models where the limit point process of the spectrum correspond to Brownian motion of the form (\ref{XXXSDE}) with other parameter values.

\subsection{Stochastic Dirac operators}

The operators $\Sineop$, $\Besselop$ and $\Huaop$ are all of the form  $\tau f=R^{-1} J f'$  where $R$ is given by $\frac{X^t X}{2\det X}$,
and $X=X_{\tch(t)}$ is a matrix-valued Brownian motion on a logarithmic time scale.

This section will show that in a certain sense an operator of this form  can be transformed into the form $J(\partial_t+\textup{``noise+drift''})$. The fact that the $\Sineb$ process can be represented as the spectrum of such a differential operator was  first conjectured in \cite{ES}.

For any invertible $2\times2$ real matrix $X$ we have $\frac{X^t X}{\det X}=J X^{-1}J^{-1} X$. Consider the new operator $\tilde \tau$ defined as
\[
\tilde \tau f(t)=X_{\tch(t)}  \tau (X_{\tch(t)}^{-1} f(t)),
\]
on the domain $\{f: X_{\tch(t)}^{-1} f(t)\in \dom(\tau)\}$. The change of variables will not change the spectrum of the operator, thus $\tilde \tau$  will have the same spectrum as $\tau$.
\emph{Heuristic} application It\^o's formula gives
\begin{align*}
\tilde \tau f&=X (2 X^{-1} J X J^{-1}) J \partial_t (X^{-1} f)=2 J X \partial_t( X^{-1} f)\\
&=2 J X(-\tch'(t)X^{-1} dX X^{-1}+\tch'(t)^2X^{-1} dX X^{-1} dX X^{-1} +X^{-1} \partial_t)f\\
&=2J(\partial_t-\tch'(t) d\cB_{\tch(t)}+\tch'(t)^2d\cB_{\tch(t)} d\cB_{\tch(t)})f.
\end{align*}
The $d\cB_{\tch(t)} d\cB_{\tch(t)}$ term simplifies to a drift term, so when $B_1$ and $B_2$ are independent we get the heuristic
$$
\tilde \tau = 2\mat{0}{-\partial_t}{\partial_t}{0}+ \frac{2}{1-t} \mat{0}{ \gamma_1 +\alpha_1 dB_1-\tfrac{\alpha_1^2}{1-t}}{0}{\gamma_2 +\alpha_2 dB_2}.
$$

\section{The soft edge operator as a canonical system}\label{s:edge}

The scaling limit of the Gaussian $\beta$-ensemble near the spectral edge $\pm 2\sqrt{n}$ is the Airy$_\beta$ point process.
This limit in the classical cases was done in \cite{TW_soft}. The general $\beta$ case was handled in \cite{RRV}.
\begin{theorem}[Soft edge limit \cite{RRV}]
Fix $\beta>0$ and let $\Lambda_n$ be a finite point process with density (\ref{betadens}). Then $n^{1/6}(2\sqrt{n}-\Lambda_n)$ converges in distribution to a  point process $\Airyb$. $\Airyb$ is the (discrete) spectrum of $\Airyop=-\partial_t^2+t+\frac{2}{\sqrt{\beta}} dB$ acting on functions $[0,\infty]\to \R$ with Dirichlet condition at 0.
\end{theorem}
Here $dB$ is white noise, and the exact definition of  the operator is given in \cite{RRV}. The fact that this operator can be represented as a self-adjoint generalized Sturm-Liouville operator on $L^2[0,\infty)$ was shown in \cite{AlexPhD}, see also \cite{Minami}.

The stochastic operator $\Airyop$ does not fit into our framework of random Dirac operators.  However the eigenvalue equation for the operator can be rewritten as a canonical system (\ref{canonical}).

Let $u_1, u_2$ be solutions of $\Airyop u_i=0$ with initial conditions $(u_i(0), u_i'(0))=(1,0)$ and $(0,1)$ for $i=1$ and $2$, respectively. This means that $u_i, u_i'$  satisfy the SDE system
\begin{align*}
du_i=u_i' dt, \qquad du_i'=u_i(\frac{2}{\sqrt{\beta}} dB+t dt).
\end{align*}
Now set $R(t)=u u^t$ where $u=(u_1(t),u_2(t))^t$. Then the solutions of the canonical system
\begin{align}\label{Airy_can}
J y'(t)=\lambda R(t) y(t), \qquad y:[0,\infty)\to \CC^2
\end{align}
are exactly of the form $y(t)=Q(t)^{-1} (v(t), v'(t))^t$ where $\Airyop v=\lambda v$ and  $$Q(t)=\mat{u_1(t)}{u_2(t)}{u_1'(t)}{u_2'(t)}.$$ See Section 8 of \cite{Remling} for more detail.

Note that the natural $L^2$ space for the canonical system (\ref{Airy_can}) is $L^2_R[0,\infty)$. Since we have $(Q^t)^{-1} R Q^{-1}=\mat{1}{0}{0}{0}$, the change of variables $y=Q^{-1} v$ will map the space $L^2_R[0,\infty)$ to the $L^2[0,\infty)$ space for the $\Airyop$ operator.

If we follow $r_\lambda(t)=y_1(t)/y_2(t)$ where $y_\lambda=(y_1,y_2)^t$ solves (\ref{Airy_can}) then we get the ODE
\begin{align}\label{softedge_car}
r'_\lambda(t)=\lambda(r_\lambda u_1(t)+u_2(t))^2.
\end{align}
The solution is strictly increasing for $\lambda>0$ and it restarts at $-\infty$ whenever it blows up to $\infty$ (with similar restarts at $-\infty$ for $\lambda<0$). The evolution of the boundary point $r_\lambda$ is similar to a hyperbolic carousel, but there is a difference in the geometry. Instead of rotating $r_\lambda$ with a fixed rate $\lambda$ along a moving center of rotation, we perform a continuously changing translation.

The evolution of the standard rate $\lambda\in \R $ translation on the boundary   is given by the solution of $r'(t)=\lambda$. Here the fixed boundary point is $\infty$. Conjugating this evolution with an isometry of $\HH$ shows that the general rate $\lambda$ infinitesimal translation acting on $\partial \HH$ is given by  $r'(t)=\lambda (a_1 r+a_2 )^2$.

\bigskip

\noindent \textbf{Acknowledgements.} We thank J.~Lagarias and  P.~Sarnak for helpful comments and references. We also thank the BIRS Research in teams program for accommodating us while working on this paper.  The first author was partially supported by the National Science Foundation CAREER award DMS-1053280. The second author was supported by the Canada Research Chair program and the NSERC Discovery Accelerator grant. 

\bibliography{szeta}

\bigskip\bigskip\bigskip\noindent
\begin{minipage}{0.49\linewidth}
Benedek Valk\'o 
\\Department of Mathematics
\\University of Wisconsin - Madison
\\Madison, WI 53706, USA
\\{\tt valko@math.wisc.edu}
\end{minipage}
\begin{minipage}{0.49\linewidth}
B\'alint Vir\'ag
\\Departments of Mathematics and Statistics
\\University of Toronto
\\Toronto ON~~M5S 2E4, Canada
\\{\tt balint@math.toronto.edu}
\end{minipage}

\end{document}